\documentclass[12pt,letterpaper]{article}
\usepackage[english]{babel}

\usepackage[style=numeric,giveninits]{biblatex}
\usepackage{csquotes}
\usepackage{wrapfig}
\usepackage{siunitx}
\usepackage{amsmath, amssymb, bm}
\usepackage{placeins}
\usepackage[T1]{fontenc}
\usepackage{subfig}
\usepackage{placeins}
\usepackage{longtable,threeparttable}
\usepackage{color}
\usepackage{float}
\usepackage{setspace}\setdisplayskipstretch{}
\usepackage{graphicx}
\usepackage{hhline}
\usepackage{multirow}
\usepackage{tabularx}
\usepackage{hyperref}
\usepackage{makecell}

\usepackage[left=1in, right=1in, top=1in, bottom=1in]{geometry}
\usepackage{setspace}

\title{Optimal Economic Operation of Liquid Petroleum Products Pipeline Systems}

\author{Elena Khlebnikova\thanks{Information Systems \& Modeling (A-1), $^\dagger$Applied Mathematics \& Plasma Physics (T-5), Los Alamos National Laboratory, Los Alamos, NM, USA. \{elenak | kaarthik | azlotnik | rbent | mewers | btasseff\}@lanl.gov. This project was carried out as part of the Advanced Grid Modeling Research program of the D.O.E. Office of Electricity. The work was conducted at Los Alamos National Laboratory under the auspices of the National Nuclear Security Administration of the U.S. Department of Energy under Contract No. 89233218CNA000001.}, Kaarthik Sundar$^*$, Anatoly Zlotnik$^\dagger$, \\ Russell Bent$^\dagger$, Mary Ewers$^*$, Byron Tasseff$^*$}

\newcommand{\edit}[1]{{\color{black}#1}}

\begin{document}

\maketitle

\begin{abstract}

The majority of overland transport needs for crude petroleum and refined petroleum products are met using pipelines. Numerous studies have developed optimization methods for design of these systems in order to minimize construction costs while meeting capacity requirements. Here, we formulate problems to optimize the operations of existing \edit{single liquid commodity} pipeline systems subject to physical flow and pump engineering constraints. The objectives are to maximize the economic value created for users of the system and to minimize operating costs. We present a general computational method for this class of \edit{continuous, non-convex nonlinear programs}, and examine the use of pump operating settings and flow \edit{allocations} as decision variables. The approach is applied to compute optimal operating regimes and perform engineering economic sensitivity analyses for a case study of a crude oil pipeline developed using publicly available data.

\end{abstract}

\section{Introduction}

\setlength\abovedisplayskip{6pt}
\setlength\belowdisplayskip{6pt}

Pipeline infrastructure is the most common mechanism used to transport hydrocarbon liquid commodities such as petroleum, crude oil, and diesel fuel \cite{9korshak2}. It is the most economically viable and the most reliable means for transporting hydrocarbons over long distances on land \cite{strogen2016environmental}. As a result, over the past century, considerable investments have been made throughout the world in continental-scale pipeline networks for transporting hydrocarbons. The sustained dependence on hydrocarbons for meeting modern demands for energy has compelled growing interest in long-standing areas of modeling and optimizing the physical and engineering function of these complex networks \cite{pharris2008argonne,rizwan2013crude,losenkov2019optimization}. Recent studies have considered both the stationary \cite{1lurie} and transient regimes \cite{2watters,4thorley2,5chaudhry}, with the predominant goals of accurately evaluating transport capacity and optimizing operational efficiency \cite{wang2012survey}.



The function of large-scale liquid transport pipelines is typically studied by characterizing physical and engineering relationships among variables such as flow rates, pressures, and applied pumping power throughout a pipeline system and its components \cite{barkhatov2017development}. The transport capacity and operating efficiency of a pipeline varies because of numerous factors arising from the condition of oil fields and processing plants (suppliers to the pipeline), hydrocarbon extraction and processing schedules, the system maintenance status, and the requirements for delivery to refineries and transport terminals (consumers of the commodity) \cite{losenokov2017optimization}. The resulting utilization of system transport capacity can vary significantly on annual, seasonal, and even daily bases. System design and operation are, therefore, conservative in order \edit{to provide flexibility to accommodate} variations in transport volumes \cite{liu2015research}.



The use of optimization for energy transport networks was pioneered in early studies on dynamic programming for gas pipeline system design \cite{wong1968optimization}. It was later developed to more efficiently and economically manage energy flow in large regional- and continental-scale electric power systems \cite{caramanis1982optimal,schweppe2013spot}, and is an integral part of wholesale electricity market operations throughout the world \cite{oneill2008towards,wood2013power,litvinov2010design}. Similar mechanisms that integrate physical and market operations for natural gas delivery have been proposed \cite{oneill1979mathematical}, and have been enabled by recent advances in computing and optimization technology \cite{rachford2000optimizing,rudkevich2017locational,zlotnik2019optimal}. The goal in these and similar studies has been to develop optimization formulations that determine continuous variables for electricity production and natural gas flow control, and to design coordination mechanisms as these systems become increasingly interdependent \cite{zlotnik2016coordinated}.


Approaches for optimal \edit{allocation of petroleum pipeline transport and refinery capacity} have been proposed over the past decades \cite{reddy2004novel,moro2004mixed,rejowski2008novel}. The majority of optimization studies on the midstream petroleum sector have focused on optimizing operating efficiency for a given flow \edit{allocation} \cite{camacho1990optimal,barreto2004optimization,cafaro2015minlp}. Contemporary studies typically seek to minimize the energy used for pumping and incorporate engineering models of varying fidelity \cite{liu2014optimal,starikov2015minimization}. A variety of studies have examined time-dependent scheduling of single commodity transport \cite{magalhaes2003crude,zhou2015dynamic,cafaro2011detailed,mostafaei2017continuous}, various network topologies (i.e. single- and multi-source) \cite{cafaro2011detailed,cafaro2015optimization}, and optimization under uncertainty \cite{mirhassani2008operational}. A recent formulation for optimizing flow schedules for crude oil pipelines has been proposed to minimize energy used for pumping and accounting for pipeline flow physics and pump engineering limitations \cite{losenkov2019optimization,losenkov2019mathematical}. The goal of this latest formulation is to adhere as closely as possible to a predetermined \edit{time-dependent} transport schedule and use storage tanks throughout the system as a controllable means to buffer the flow regime. 



In this study, we build on recent, previous work and develop a general optimization-based methodology for synthesizing \edit{optimal steady-state flow allocation solutions} for branched \edit{single liquid commodity} pipeline systems based on engineering economic criteria. We account for the behavior of participants in the petroleum and petroleum products supply chain \cite{soud14agent} using a model of a two-sided single auction market, similar to formulations that were previously proposed for natural gas pipeline systems \cite{rudkevich2017locational,zlotnik2019optimal}. Specifically, we formulate a problem for determining a steady-state \edit{flow allocation} and pump operating settings for a single liquid commodity product pipeline where there is no a priori nominal transport allocation, and the injections and withdrawals are decision variables. We assume that each supplier or consumer that is a user of transportation services provided by the pipeline expresses \edit{its economic} position as a range of possible (feasible) flow rates (either \edit{for injection} or withdrawal) at a given price level. These positions are interpreted as bids into an optimization-based single auction market whose solution determines the optimal flow configuration of the system. \edit{We restrict our modeling and optimization goals to a fungible transportation mode, in which we assume that product shipments are indistinguishable, or fungible, and consumers are open to receive any product coming from any source.} The problem is solved subject to equality constraints that represent the physics of weakly compressible fluid flow through pipes and a model of the mechanics of centrifugal pumping equipment that is used to propel flow through the system. We develop a simplified basic model of pumping station performance based on contemporary models of variable frequency drive machinery \cite{14Samolenkov,zemenkov2017}, as well as flow equations that account for anti-turbulent additives that have now become standard in the petroleum pipeline sector \cite{7Leibenzon_eq}. The problem incorporates inequality constraints that reflect the operating limits on flows, pressures, and other state and control variables.
When the objective function of the continuous optimization problem maximizes the payments by consumers minus the receipts of suppliers and minimizes the cost of energy used in pump operation, \edit{we propose that} the optimal solution ensures that the pipeline operates in the feasible configuration that provides the greatest economic value to users of the system. Our goal here is to pose the problem and examine the outcome of the optimization on a synthesized yet realistic case study. Additionally, we propose that the Lagrange multipliers that are obtained as output of the optimization solution are interpretable as values of the transported commodity throughout the network. We perform a perturbation analysis on the economic parameters in order to understand the sensitivity of these locational prices to changes in inputs to the mathematical problem. Ultimately, our engineering economics approach could be used by liquid pipeline managers as a decision support tool for integrated financial and physical operations of their systems to be responsive to customer requirements and energy prices. The application of efficient and responsive liquids pipeline scheduling could be used to mitigate the carbon dioxide emissions of crude oil transportation and refining \cite{jing2019global}.

The remainder of the paper is organized as follows. In Section \ref{sec:flow}, we review modeling equations for weakly compressible fluid flow in a pipe. Then in Section \ref{sec:pipeline}, we review network modeling for pipelines, and integrate the flow equations with basic engineering modeling of pumping machinery. Next, Section \ref{sec:formulation} is devoted to the formulation of the optimization problem to address the stated goals and requirements. Then, Section \ref{sec:results} provides details on the computational implementation of our approach, and presents results for a test study synthesized based on an actual pipeline system and prevailing economic data, as well as a sensitivity analysis. Comprehensive tables that specify the system model, as well as the inputs and outputs for the optimization problem, are provided. We conclude with a discussion of the performance of our method, follow-on work, and potential applications. 
\section{Flow Equations} \label{sec:flow}


We now review the dynamic flow equations for weakly compressible fluid flow in long pipes and the assumptions made about the flow regime. We also explain the reduction to steady-state flow equations which relate head decrease to flow rate in the presence of drag-reducing additives.

\subsection{Unsteady Fluid Flow in a Pipeline}
\label{subsec:unsteadyflow}
Pipeline fluid flow in the transient regime, in which hydraulic parameters such as flow rate, pressure, speed, and temperature depend not only on the space coordinate, $x$, along the axis of the pipeline, but also on time $t$, is commonly referred to as unsteady flow. Unsteady processes in pipelines can occur, for example, during initiation and shutdown of transport operations of a main pipeline, or in the case of pipeline rupture that results in high dynamic transients \cite{hajossy2014cooling}. The model of unsteady flow of a weakly compressible fluid such as petroleum or a liquid petroleum product in a pipeline is described in several studies \cite{4thorley2,wylie1993fluid}. We describe the propagation of a non-stationary process along the length of a pipeline using the following equations for mass and momentum conservation \cite{1lurie}:
\begin{subequations} 
\begin{align}
\frac{\partial p(t,x)}{\partial t} + c^{2} \rho \frac{ \partial u(t,x)}{\partial x} &= 0, \label{eq:pde0a}\\ 
\rho\frac{\partial u(t,x)}{\partial t} + \frac{\partial p(t,x)}{\partial x} & = -\lambda(Re,\epsilon) \rho \frac{ (u(t,x))^2}{2 D} -  g \sin(\alpha(x)) \rho.  \label{eq:pde0b}
\end{align}
\end{subequations}
The variables $p(t,x)$ and $u(t,x)$ in equations \eqref{eq:pde0a}-\eqref{eq:pde0b} model the pressure and flow velocity as functions of time and space. Because we consider a weakly compressible fluid, we assume that the fluid density $\rho$ is constant and uniform throughout the pipe and network model. The parameters include the hydraulic resistance coefficient $\lambda(Re,\epsilon)$, which depends on the Reynolds number $Re = uD/\nu$, the flow velocity $u = u(x, t)$, and the pipe cross-sectional area $A$. Here $D$ is the pipeline diameter and $\epsilon = \Delta/D$ is the relative roughness of the pipe inner surface. The spatially distributed parameter $\alpha(x)=dz/dx$ is the inclination angle of the pipe axis with respect to the horizontal, where $z(x)$ is elevation of the pipe axis. The \edit{coefficient $c$} is the speed of wave propagation in the pipe, which is constant given our assumption of constant density, and $g$ is gravitational acceleration. To determine the value of the hydraulic resistance coefficient $\lambda$ for different flow regimes, we use the equations of Stokes, Blasius, Altshul, and Shifrinson \cite{6lipovka},
\begin{align}\label{eq:friction}
\lambda & =\begin{cases}
\frac{64}{Re},  Re<2320\\
\frac{0.3164}{{Re}^{0.25}},  2320<Re<{10}^{5}\\
0.11{(\frac{68}{Re}+\frac{\Delta}{D})}^{0.25},    {10}^{5}<Re<500\frac{\Delta}{D}\\
0.11{(\frac{\Delta}{D})}^{0.25},   Re>500\frac{\Delta}{D}.
\end{cases}
\end{align}

The solution to the system of equations \eqref{eq:pde0a}-\eqref{eq:pde0b} requires specification of initial conditions that characterize the distributions of pressure $p(0,x)$ and flow velocity $u(0,x)$ at the initial time $t = 0$, boundary conditions that reflect the setting at the endpoints $x = 0$ and $x = L$ of a pipe, as well as the conjugation conditions. The latter conditions are \edit{imposed} in some intermediate sections $x = x_0$ and reflect the effect of the operation of pipeline valves, bends for the selection or pumping of liquid, and other events that occur at such locations.

\subsection{Steady Fluid Flow in a Pipeline} \label{subsec:steadyflow}

In the steady state, the pressure and flow velocity variables in each section of pipeline remain constant as functions of time, and are related according to the mechanical energy conservation equation
\begin{align}\label{eq:steadyf}
\frac{d p}{\rho} + \lambda \frac{dx}{D}\frac{u^2}{2} + d \frac{u^2}{2} + gdz = 0.
\end{align}

Because the fluid is weakly compressible and homogeneous, $\rho \equiv $ constant and hence $d p/d x = 0$. Since the pipe diameter is constant, the velocity $u(x)$ of fluid remains constant throughout the pipe. In equation \eqref{eq:steadyf}, the value of $d p/\rho$ represents the work of a unit of mass of fluid moving along the area $dx$. This work is spent on overcoming the resistive frictional forces of turbulence, the change in fluid kinetic energy, and on lifting fluid over a height difference of $\Delta z$. After integrating the equation of mechanical energy conservation between sections $x = x_1$ and $x = x_2$ of a pipe, where ${L}_{1-2}$ is the distance between sections $x_1$ and $x_2$, the relation between endpoint pressure variables $p_1\equiv p(x_1)$ and $p_2\equiv p(x_2)$ and the flow velocity $u$ is given by
\begin{align}\label{eq:bernoulli0}
\frac{p_1 - p_2}{\rho g} &= \lambda \frac{{L}_{1-2}}{D}\frac{{u}^{2}}{2g}+\Delta z.
\end{align}

The values $H_1\equiv p_1/(\rho g)$ and $H_2\equiv p_2/(\rho g)$ represent the heights the liquid is lifted to under the pressures $p_1$ and $p_2$ at locations $x=x_1$ and $x=x_2$, respectively. This value is referred to as head. The value of $\Delta H \equiv (p_1 - p_2)/(\rho g)$ is called total head loss. In the general case, the causes of total head loss are local resistances and fluid speed changes. The first term on the right hand side of equation \eqref{eq:bernoulli0} is the Darcy-Weisbach formula  \cite{brown2003history}.
A general equation that is used to determine the coefficient of hydraulic resistance $\lambda$ is given as follows:
\begin{align}\label{eq:lambda}
\lambda &= \frac{A}{Re^m},
\end{align}
where cross-sectional area $A$ \edit{and} a value corresponding to the fluid motion regime $m$ are constant values. Substituting equation \eqref{eq:lambda} into the Darcy–Weisbach equation and taking into account the relation $Re = 4Q/(\pi D \nu)$, where $Q$ is the volumetric flow rate and $\nu$ is the petroleum kinematic viscosity, we obtain the Leibenzon formula \cite{7Leibenzon_eq}:
\begin{align}\label{eq:leibenzon}
h_{fr}&=1.02   \beta \frac{Q^{2-m} \nu^m} {D^{5-m}}.
\end{align}
For turbulent flow in smooth pipes, we use standard parameters of $m$ = 0.25 and $\beta$ = 0.0246 $\textnormal{s}^2/\textnormal{m}$ \cite{7Leibenzon_eq}.
The Leibenzon formula is widely used in cases where the dependence of $h_{fr}$ on $Q$ is expressed explicitly. The described mathematical model in \cite{7Leibenzon_eq} enables the prediction of a sufficiently accurate hydraulic resistance resulting from the viscous friction of turbulent flow in a pipeline. This model is based on the linear transformation law of the coefficients in the generalized Leibenzon formula. The conversion is a function of an input parameter that reflects the degree of laminarization of turbulent flow in pipes with smooth walls. In this manner, the head loss caused by pipe flow can be represented without \edit{dependence on the Reynolds} number and, furthermore, the parameter $m$ can be adjusted to reflect the laminarization effects of anti-turbulent drag reduction additives that are widely used to facilitate petroleum transport.

\section{Liquid Pipeline Network Modeling} \label{sec:pipeline}
In this study, we present a basic and generalizable representation of a large-scale liquid pipeline network, which is based on our previous work \cite{khlebnikova2020optimization}. The flow of fluid through each type of component in the proposed model is represented with physical laws and engineering properties. These physical constraints, together with an objective function and inequality constraints that restrict the variable values, define an optimization problem that is described in Section \ref{sec:formulation}. The formulation is developed to be solvable using general purpose optimization software and to be applicable to a variety of network topologies and components. Our intention is to optimize the function of very large systems on a short (operational) time scale, and the modeling includes the basic elements of pipes, junctions, and pumps. In this study, we do not consider integer variables for activation of components or changes in system topology. We assume that component subsystems such as pumping station controls are managed to track the overall system setpoints.
\subsection{Network Modeling}
In this model, a liquid's transport pipeline consists of pipe sections that are connected at junctions, through which flow of a weakly compressible fluid is propelled by centrifugal pumping equipment. The representation includes a model of the static network system (which does not change), state and control variables for the system, and all component and fluid parameters. Our modeling facilitates optimization of flow \edit{allocations} over large scales, which results in a nonconvex, nonlinear problem that is formulated in Section \ref{sec:formulation}.
We use several assumptions in formulating the network model. Generally, we assume that the system is designed to operate continuously in steady state for balanced flows and that the flow direction on each pipeline segment is determined beforehand. Thus, we do not model transient processes (such as those arising from cavitation) and assume that any transients stabilize quickly after control changes are made concurrently throughout the system. The movement of liquid is modeled as one-dimensional turbulent flow. Further, we assume that the liquid commodity is homogeneous, so that properties such as density and viscosity are constant and uniform throughout the system. Following precedent, we suppose that the inaccuracies caused by these homogeneity assumptions are minor  \cite{losenkov2019optimization}. We do not account for any changes in temperature and do not model maintenance scenarios, like pigging in the system.
We represent the pipeline network mathematically as a set of junctions (or nodes) indexed as $j \in \mathcal N$, where $\mathcal N_S \subset \mathcal N$ and $\mathcal N_D \subset \mathcal N$ represent subsets of the nodes with supplies and consumptions of the commodity, respectively. The model also includes pipes (or edges) $(i,j) \in \mathcal E$ that connect pairs of nodes in $\mathcal N$, where $\mathcal E_E \subset \mathcal E$ and $\mathcal E_P \subset \mathcal E$ denote subsets of the edges that represent pipes and pump stations, respectively. We assume that the pair $(\mathcal N,\mathcal E)$ defines a connected graph. We use the notation $N=|\mathcal N|$ for the number of nodes, $E=|\mathcal E_E|$ for the number of pipes, and $P=|\mathcal E_P|$ for the number of pumps. Next, we present the equalities that represent physical and engineering relationships for \edit{the} flow through pipeline sections, junctions, and pumps, respectively, which are modeled as equality constraints in the optimization formulation presented thereafter.
\subsection{Pipe Flow Model}
The general purpose of liquid pipeline hydraulic calculations is to determine the head losses of moving fluid through the pipeline. In general, the total head loss between the ends of a pipe is equal to the sum of losses caused by friction and the head difference arising from pipeline elevation. In this study, we use the Bernoulli equation,
\begin{align}\label{eq:bernoulli1}
H_i-H_j&=z_j-z_i+1.02\beta_{ij} \frac{Q_{ij}^{2-m} \nu^m} {D_{ij}^{5-m}} L_{ij}, \quad \forall (i,j) \in \mathcal E_E,   
\end{align} 
where the head losses are defined using the Leibenzon equation \eqref{eq:leibenzon}. Here, $H_j$ and $z_j$ denote the head and elevation at a node $j \in \mathcal N$, and $Q_{ij}$, $L_{ij}$, and $\beta_{ij}$ model the flow rate through a pipe $(i,j) \in \mathcal E_E$, and its length and friction coefficient, respectively. 
\subsection{Nodal Balance Model}
As standard in the modeling of network flow systems, we include a Kirchhoff-like flow balance law. We use shorthand to denote the sets of nodes from which flows are incoming to, and to which flows are outgoing \edit{from} a node $j \in \mathcal N$ by $\partial_{+}j = \{i \in N | (i,j) \in  \mathcal E\}$ and $\partial_{-}j  = \{k  \in N | (j, k) \in \mathcal E\}$. For the liquid pipeline systems studied here, nodal volume flow balance is given as
\begin{align}\label{eq:flowbal}
\displaystyle \sum_{i \in \partial_{+}j} Q_{ij} - \displaystyle \sum_{k \in \partial_{-}j} Q_{jk} = s_j - d_j, \, \quad \forall j \in \mathcal N,
\end{align}
where $s_j$ and $d_j$ denote injection and withdrawal of the commodity into the pipeline at the node $j \in N$. Although these values are in general not concurrently nonnegative, we use these two distinct variables to promote generalizability of optimization formulations for network flow problems that are interpreted as auction markets for a commodity and/or transport capacity.

\subsection{Pump Model} \label{subsec:pump}
In liquid pipeline transportation, the flow of the commodity through the system is actuated by pumping machinery, which is located at pumping stations throughout the system. The major actions that can be taken to conserve energy used for the operation of legacy systems are the use of adjustable variable frequency drives at pumping stations, the use of smooth pipe coatings, anti-turbulent additives such as drag reduction agents mixed with the transported product, replacement and reconstruction of pumping stations and their drives, the use of static capacitors, and the expansion of automation facilities and management capabilities. Although pumping machinery is typically located at large pumping stations with multiple units and complex topology, our goal is to optimize flows on overall system-wide scales. We assume simple relationships between a minimal set of physical and mechanical variables that represent pump action and consider pump stations to consist of several standard pump units arranged in series. We now review the engineering assumptions and  relations that we synthesize into the variable pump drive modeling equations, which are subsequently included as constraints in  optimization formulations.
The essence of frequency control is the use of frequency regulators in adaptive electric motor drives. Such components enable efficient control of an electric pump depending on operating conditions. The effect is achieved by changing the frequency and amplitude of the three-phase voltage supplied to the electric motor. Here we describe the relationship between rotation frequency, pump efficiency, and flow. We drop the $(i,j)$ indices for ease of exposition. The dependence of pump efficiency on flow at nominal rotational speed is modeled by a pump-specific characteristic curve $\eta = f(Q)$. This curve is often approximated analytically \cite{2012Vyazunov, kocher2016analytical} using the quadratic function 
\begin{align}\label{eq:pumpeff0}
\eta^{pump} &= b_0+b_1 Q+b_2 Q^2,
\end{align}
where $b_0$, $b_1$, and $b_2$ are approximation coefficients obtained by the least squares method and $Q$ is the flow rate through the pump in $m^3/h$.
The change of centrifugal pump efficiency at variable speed was studied in \cite{2012Vyazunov, georgescu2014estimation, marchi2012assessing, 12Grishin}. Affinity Laws can be used to define the pump efficiency. The values of the actual and nominal flow rates $Q$ and $Q^{nom}$ and heads $H$ and $H^{nom}$ when operating at actual and nominal rotational pump speeds $\omega$ and $\omega^{nom}$, respectively, are related by the ratios
\begin{subequations}
\begin{align}
\frac{Q}{Q^{nom}} &= \frac{\omega}{\omega^{nom}}, \label{eq:pumpratio1} \\
\frac{H}{H^{nom}} &= \left( \frac{\omega}{\omega^{nom}}\right)^2. \label{eq:pumpratio2}
\end{align}
\end{subequations}
Thus, we obtain equation \eqref{eq:pumpeff1} for the pump efficiency as a function of rotational speed:
\begin{align}\label{eq:pumpeff1}
\eta^{pump} &= b_0+b_1\frac{Q^{nom}}{\omega^{nom}}\omega+b_2\left( \frac{Q^{nom}}{\omega^{nom}}\right)^2\omega^2
\end{align}
The laws of affinity require the use of a nominal point on the curve to determine the new rotational speed. Most references do not provide a clear description of the method that determines the base point, and the wrong choice of this point can significantly affect the results. The following equation, given by \cite{sarbu1998energetic}, provides an analytical relationship between actual and nominal rotational speeds $\omega$ and $\omega^{nom}$ and corresponding efficiencies $\eta$ and $\eta^{nom}$:
\begin{align}\label{eq:nominaleff}
\eta &= 1 - (1 - \eta^{nom})\left(\frac{\omega^{nom}}{\omega}\right)^{0.1}.
\end{align}
Therefore, for large pumps the change in efficiency can be neglected if the frequency is within 33\% of its nominal value \cite{sarbu1998energetic}. It follows that this approximation is valid when the points of the pump operating curve will maintain the same efficiency. However, this does not mean the pump will perform at the same efficiency when installed in the pipeline system, because the operating point is determined by the intersection of the pump curve with the system curve.
Grishin \cite{12Grishin} proposed to model the dependence of pump efficiency on flow at a constant speed using a second-degree polynomial:
\begin{align}\label{eq:pumpeff2}
\eta^{pump} & = b_1 Q^{nom}-b_2 (Q^{nom})^2,
\end{align}
which can be rewritten as
\begin{align}\label{eq:pumpeff3}
\eta^{pump} & = 2\frac{\eta^{nom}}{Q^{nom}}Q-\frac{\eta^{nom}}{(Q^{nom})^2}Q^2,
\end{align}
where $\eta^{nom} = b_1^2/(4b_2)$ and $Q^{nom} = b_1/(2b_2)$. It follows that the dependence of efficiency on the actual flow rate can be approximated with
\begin{align}\label{eq:pumpeff4}
\eta^{pump} & = 2\frac{\eta^{nom}}{Q^{nom}}Q \frac{\omega^{nom}}{\omega}-\frac{\eta^{nom}}{(Q^{nom})^2}Q^2 \left(\frac{\omega^{nom}}{\omega}\right)^2.
\end{align}
Equation \eqref{eq:pumpeff4} can then be transformed using the technique of completing the square to yield
\begin{align}\label{eq:pumpeff5}
\eta^{pump} & = \eta^{nom} - \left(Q - Q^{nom}\frac{\omega}{\omega^{nom}}\right)^2 \cdot \frac{\eta^{nom}}{(Q^{nom})^2} \cdot \left(\frac{\omega^{nom}}{\omega}\right)^2.
\end{align}
Thus, by modulating the power applied to a pump $(i,j)\in \mathcal E_{p}$, it is possible to make the engine drive shaft rotational speed $\omega_{ij}$ lower or higher than a nominal value $\omega_{ij}^{nom}$. The relationship between the drive frequency, head difference, and through flow for a variable frequency drive pump \cite{10shabanov1,12Grishin} can then be approximated by equations
\begin{subequations}
\begin{align}
H_{j} - H_i  = a_{ij}^0 \left(\frac{\omega_{ij}}{\omega_{ij}^{nom}}\right)^2 - a_{ij}^1 Q_{ij}^2, & \quad \, \forall (i,j) \in  \mathcal E_P, \label{eq:pumpflow}\\
\eta_{ij}  = \eta_{ij}^{nom} - \left(\frac{Q_{ij}}{Q_{ij}^{nom}} - \frac{\omega_{ij}}{\omega_{ij}^{nom}}\right)^2 \eta_{ij}^{nom} \left(\frac{\omega_{ij}^{nom}}{\omega_{ij}}\right)^2 , & \quad \, \forall (i,j) \in  \mathcal E_P, \label{eq:pumpefficiency}
\end{align}
\end{subequations}
where $a_{ij}^0$ and $a_{ij}^1$ are constant parameters and  $Q_{ij}^{nom}$, $\omega_{ij}^{nom}$, and $\eta_{ij}^{nom}$ are nominal values of flow rate, drive shaft frequency, and pumping efficiency, respectively, for pump $(i,j)\in \mathcal{E}_{p}$. 
Equations \eqref{eq:pumpflow} and \eqref{eq:pumpefficiency} provide a minimal relationship between the input and output head and through flow rate of an electric drive pump and its drive shaft frequency and efficiency. The energy used for operating the pump $(i,j)\in \mathcal{E}_{p}$ in a given configuration can be approximated by \cite{11Shabanov2}
\begin{align}\label{eq:pumpenergy}
E_{ij} &=  \rho g Q_{ij} \frac{H_{j}-H_{i}}{\eta_{ij} \eta_{em} \eta_{mt}}, \quad \, \forall (i,j) \in  \mathcal E_P,
\end{align}
where $\rho$ is the liquid's density, and $\eta_{em}$ and $\eta_{mt}$ are the electric motor and mechanical transmission efficiency of the drive, respectively, which are both assumed to be constant and uniform. We use the same values of these constant efficiencies for pumps throughout the system. Finally, it is assumed that the pump is operating in its designed forward direction so that $H_j>H_i$.

\section{Optimization Formulation} \label{sec:formulation}

We now present a mathematical optimization formalism that reflects the business goals of a liquid pipeline system manager and ensures operational and safety requirements. We assume that these goals are to 1) maximize the profit gained from providing a transportation service for a liquid commodity that is shipped from sellers to buyers on a pipeline network; 2) minimize the cost of operating the pipeline system; and 3) to optimize the previous criteria while strictly maintaining the values of physical and mechanical variables within allowable engineering limits. To do so, we formulate several optimization problems that consist of 1) objective functions to reflect the optimized quantities; subject to 2) the equality constraints \eqref{eq:leibenzon}-\eqref{eq:flowbal} and \eqref{eq:pumpflow}-\eqref{eq:pumpefficiency} that describe steady-state pipe flow physics and pump drive action; and 3) inequality constraints that reflect limits on pump operations and physical variables. Here we examine several different formulations to evaluate the economic benefit of optimizing the commodity flow \edit{allocation}.

\subsection{Objective Function}
The objective function is a difference of two cost terms that give the total economic value provided by the pipeline in a given time period. The first and most significant term \edit{is the economic objective}, $J_E$, which defines the economic value of transportation that is provided to users of the pipeline system and is given as the difference between the payments for withdrawals and receipts for injections,
\begin{align}\label{eq:costobj}
J_{E} &=\sum_{j\in \mathcal N_D} t_j d_j - \sum_{j\in \mathcal N_S} r_j s_j. 
\end{align}

In equation \eqref{eq:costobj}, $s_j$ and $d_j$ are supplied injections and demanded withdrawals, respectively, and $r_j$ and $t_j$ are prices offered by sellers and buyers of the commodity at each node $j \in \mathcal N$. The second term defines the cost of operating pumps on the system. We assume that the pumps are operating continuously, so that maintenance costs are fixed and the primary operating cost is from the electricity used by pump drives. The pumping energy is defined in terms of the relative pump shaft rotation speed, i.e., the ratio of the actual speed to the nominal speed. From \eqref{eq:pumpenergy}, the form of the objective function value for the cost of energy expended for pumping system-wide is
\begin{align}\label{eq:energyobj}
J_{O} &= \rho g  \sum_{(i,j)\in\mathcal E_P}\frac{H_{j}-H_{i}}{\eta_{ij} \eta_{em} \eta_{mt}} Q_{ij} c_{ij},
\end{align}
where $c_i$ is the price of electric energy at the pump station header node. The total objective function that we optimize is then given as the difference of the \edit{economic objective \eqref{eq:costobj} and operating objective \eqref{eq:energyobj}}, which gives the \edit{total} instantaneous economic value (over a given time period) produced by the pipeline system,
\begin{align}\label{eq:combinedobj}
J_{P} &=  J_E - J_O.
\end{align}
\subsection{Inequality Constraints}
\setlength\abovedisplayskip{6pt}
\setlength\belowdisplayskip{6pt}
Changes in injection and withdrawal of liquid from a pipeline can cause pressures to exceed maximum and minimum allowable values, creating ruptures and cavitation. Analysis of deviations in operating modes of a pipeline system is also important for preventing damage to or premature wear on pumping machinery and preventing inefficient electricity consumption. Therefore, inequality constraints that limit the flow, head, and pump operating points are needed to reflect operating requirements for the system. First, feasibility of the physical flow and head solution is enforced by delimiting the fluid flow rate on each pipe and the allowable head at each junction:

\vspace{-4ex}
\singlespacing
\begin{subequations}
\begin{align}
Q_{ij}^{min} \leq Q_{ij} \leq Q_{ij}^{max}, & \quad \forall \,  (i,j) \in  \mathcal E_E, \label{eq:pipeflowlim} \\
H_{j}^{min}  \leq H_{j} \leq H_{j}^{max}, & \quad \forall \, j \in  \mathcal N. \label{eq:nodeheadlim}
\end{align}
\end{subequations}

In addition, the acceptable operating states of variable drive electric pumps are defined by inequalities for effective flow, rotation speed, pumping efficiency, and head difference between inlet and outlet:

\vspace{-4ex}
\singlespacing
\begin{subequations}
\begin{align}
0.8Q_{ij}^{nom} \leq Q_{ij} \leq 1.2Q_{ij}^{nom}, & \quad \forall \,  (i,j) \in  \mathcal E_{p}, \label{eq:pumpflowlim} \\
0.7\eta_{ij}^{nom} \leq \eta_{ij} \leq \eta_{ij}^{nom},& \quad \forall \, (i,j) \in  \mathcal E_{p}, \label{eq:pumpefflim} \\
0.8 \omega_{ij}^{nom} \leq \omega_{ij}\leq 1.2 \omega_{ij}^{nom}, & \quad \forall \, (i,j) \in  \mathcal E_{p}, \label{eq:pumpfreqlim} \\
H_{ij}^{min} \leq H_{j}-H_{i} \leq H_{ij}^{max}, & \quad \forall \, (i,j) \in  \mathcal E_{p}. \label{eq:pumpheaddifflim}
\end{align}
\end{subequations}

Furthermore, if the flow configuration of the system is optimized, then the injections and withdrawals are decision variables and must be constrained with 

\vspace{-4ex}
\singlespacing
\begin{subequations}
\begin{align}
s_{j}^{min} \leq s_{j} \leq s_{j}^{max}, & \quad  \, \forall j \in  \mathcal N_S, \label{eq:supplylim} \\
d_{j}^{min} \leq d_{j} \leq d_{j}^{max}, & \quad \, \forall j \in  \mathcal N_D. \label{eq:demandlim}
\end{align}
\end{subequations}

\subsection{Optimization Formulations} \label{subsec:formulations}
We next present three formulations for optimizing pipeline operations. The formulations extremize the three objective functions in equations \eqref{eq:costobj}, \eqref{eq:energyobj}, and \eqref{eq:combinedobj} subject to nonlinear equality constraints and box constraints on the state and control variables that guarantee feasibility. The first formulation is similar to the optimization goals that have been examined in previous work \cite{losenkov2019optimization, camacho1990optimal}, where the flow schedule \edit{or allocation} is predetermined so that the quantities $s_j$ and $d_j$ for node $j \in \mathcal N_S \cup \mathcal N_D$ are fixed and the objective is to minimize the (operational) cost of pumping energy.

\singlespacing
\noindent Formulation 1:
\begin{longtable}{ l p{4in} } 
min & $J_O$ = cost of energy used for pump operation \eqref{eq:energyobj} \\
s.t. & Pipe flow equation \eqref{eq:bernoulli1}\\
& Nodal flow balance with fixed nominations \eqref{eq:flowbal}\\
& Pump relations \eqref{eq:pumpflow}-\eqref{eq:pumpefficiency}\\
& Flow and head feasibility \eqref{eq:pipeflowlim}-\eqref{eq:nodeheadlim}\\
& Pump configuration feasibility \eqref{eq:pumpflowlim}-\eqref{eq:pumpheaddifflim}.
\end{longtable}

\noindent In the second formulation, the supplies and deliveries of the commodity $s_j$ and $d_j$ for nodes $j \in \mathcal N_S \cup \mathcal N_D$ are decision variables, so that the flow \edit{allocation} is optimized as well. The objective is to maximize the economic surplus created by transporting the commodity from points of low price to high price and also maximize throughput, or capacity use.

\singlespacing
\noindent Formulation 2:
\begin{longtable}{ l p{4in} } 
max & $J_E$ = economic value of commodity transport \eqref{eq:energyobj} \\
s.t.& Pipe flow equation \eqref{eq:bernoulli1}\\
& Nodal flow balance \eqref{eq:flowbal}\\
& Pump relations \eqref{eq:pumpflow}-\eqref{eq:pumpefficiency}\\
& Flow and head feasibility \eqref{eq:pipeflowlim}-\eqref{eq:nodeheadlim}\\
& Pump configuration feasibility \eqref{eq:pumpflowlim}-\eqref{eq:pumpheaddifflim} \\
& Shipper nomination limits \eqref{eq:supplylim}-\eqref{eq:demandlim}.
\end{longtable}

\noindent Finally, we synthesize the first two formulations so the flow \edit{allocation} is adjusted in order for the total volume of transported commodity to increase when the economic value that this brings to users of the system is greater than the incremental cost of pumping, subject to feasibility constraints. Crucially, all factors are considered system-wide, as the relationship between the physical flow equations, pump mechanics, and energy expended is highly nonlinear. 

\singlespacing
\noindent Formulation 3:
\begin{longtable}{ l p{4in} } 
max	&  $J_P$ = total economic value of pipeline operation \eqref{eq:combinedobj}\\
s.t.& Pipe flow equation \eqref{eq:bernoulli1}\\
& Nodal flow balance \eqref{eq:flowbal}\\
& Pump relations \eqref{eq:pumpflow}-\eqref{eq:pumpefficiency}\\
& Flow and head feasibility \eqref{eq:pipeflowlim}-\eqref{eq:nodeheadlim}\\
& Pump configuration feasibility \eqref{eq:pumpflowlim}-\eqref{eq:pumpheaddifflim} \\
& Shipper nomination limits \eqref{eq:supplylim}-\eqref{eq:demandlim}.
\end{longtable}

\noindent  Mathematically, these formulations are expressed as nonlinear programs of the form 

\vspace{-4ex}
\singlespacing
\begin{align}
\label{eq:ipoptform}
\begin{array}{ll}
\displaystyle\max_{x \in \mathbb{R}^n} &f(x) \\
s.t. & g^L \leq g(x) \leq g^U\\
& x^L \leq x \leq x^U\\
\end{array}
\end{align}

In formulation \eqref{eq:ipoptform} $x \in \mathbb{R}^n$ are the optimization variables (possibly with lower and upper bounds, $x^L\in(\mathbb{R}\,\cup\,\{-\infty\})^n$ and $x^U\in(\mathbb{R}\,\cup\,\{+\infty\})^n$); $f:\mathbb{R}^n\longrightarrow\mathbb{R}$ is the objective function, and $g: \mathbb{R}^n \longrightarrow\mathbb{R}^m$ defines constraints. The functions $f(x)$ and $g(x)$ can be linear or nonlinear and convex or non-convex, and should be twice continuously differentiable. The constraints $g(x)$ have lower and upper bounds $g^L\in (\mathbb{R}\cup\{-\infty\})^m$ and $g^U\in(\mathbb{R}\,\cup\,\{+\infty\})^m$. The equality constraints of the form $g_i(x)=\overline{g}_i(x)$ can be specified by setting $g_i^L=g_i^U=\overline{g}_i$. The constraint functions $g$ can then be defined to reflect the equalities \eqref{eq:bernoulli1}, \eqref{eq:flowbal}, and \eqref{eq:pumpflow}-\eqref{eq:pumpefficiency}, which are unalterable physical constraints, and the pump head limits \eqref{eq:pumpheaddifflim} as inequalities. The equations \eqref{eq:pipeflowlim}-\eqref{eq:nodeheadlim}, \eqref{eq:pumpflowlim}-\eqref{eq:pumpfreqlim}, and \eqref{eq:supplylim}-\eqref{eq:demandlim} can be expressed as box constraints on the variables. 

Our goal is to integrate the business functions and physical operations of the system. In the above formulations, we make a distinction between state and control variables. Whereas state variables reflect the physical state of the system, they cannot be directly controlled. The control variables are those that can be directly manipulated by the pipeline operators. The state variables are the volumetric flow rate $Q_{ij}$ on each pipe $(i,j)\in \mathcal E_E$ in m$^3/$s; the head $H_j$ at each node $j\in\mathcal N$ in m; and the actual efficiency $\eta_{ij}$ of each pump $(i,j)\in\mathcal E_P$. The control variables are the drive rotation speed $\omega_{ij}$ of each pump $(i,j)\in \mathcal E_P$ in Hz; the optimized demand $d_j$ at each node $j \in \mathcal N_D$ in m$^3/$s; and the optimized petroleum supply $s_j$ at each node $j\in\mathcal N_S$ in m$^3/$s. Note that the last two variables are included only in Formulations 2 and 3, whereas in Formulation 1 they are fixed parameters.

The parameters that determine the physical properties and limitations of the pipeline system are the resistance coefficient $\beta_{ij}$ in the Leibenzon equation, and length $L_{ij}$ and diameter $D_{ij}$ in m of each pipe $(i,j)\in \mathcal E_E$; pipe elevation $z_j$ at each node $j$ in m; the petroleum density $\rho$ and viscosity $\nu$; the pump coefficients $a_{ij}^0$ and $a_{ij}^1$ and nominal rotation speed $\omega_{ij}^{nom}$, nominal flow rate $Q_{i,j}^{nom}$, and nominal efficiency $\eta_{ij}^{nom}$ of each pump $(i,j)\in \mathcal E_P$; the minimum and maximum flow rates $Q_{ij}^{min}$ and $Q_{ij}^{max}$ on each pipe $(i,j)\in \mathcal E_E$; the minimum and maximum head $H_j^{min}$, $H_j^{max}$ on each node $j\in \mathcal N$; the minimum and maximum head difference $H_{ij}^{min}$ and $H_{ij}^{max}$ between discharge and suction of each pump $(i,j)\in \mathcal E_P$; and the price of electricity $c_j$ at pump header nodes $i$ for $(i,j)\in\mathcal E_P$.

Formulations 2 and 3 contain additional parameters for the prices offered by users of the system and the bound values for the inequalities \eqref{eq:supplylim} and \eqref{eq:demandlim} on shipper flow rates. These are the prices $t_j$ and $r_j$ offered by the seller (resp. buyer) at node $j$ in \$$/$m$^3$; the minimum and maximum flow rates $s_j^{min}$ and $s_j^{max}$ in m$^3/$h at which the seller offers the commodity at the price $t_j$ for each node $j\in\mathcal N_S$; and the minimum and maximum flow rates $d_j^{min}$ and $d_j^{max}$ in m$^3/$h at which the buyer is willing to take the commodity at the price $r_j$. We wish to emphasize that when Formulations 2 and 3 are intended to represent auction markets for liquid products transport, the triples ($t_j,s_j^{min}, s_j^{max}$) and ($r_j,d_j^{min},d_j^{max}$) can be interpreted as bids to sell (resp. buy) the commodity at the given rates.

\subsection{Dual Variables for Marginal Pricing}

A major advance in the use of optimization to manage energy transport networks was made with the development of engineering economics for electric power system markets \cite{caramanis1982optimal,schweppe2013spot}. Integrated decision making to determine generation setpoints and electricity spot pricing based on formal optimization is now at the core of wholesale electricity market operations throughout the world \cite{oneill2008towards,wood2013power,litvinov2010design}. The application of optimization for this purpose is based on a key enabling insight. When the objective function of the resource allocation problem reflects the economic value that is provided by the energy delivery network, then the dual variables (Lagrange multipliers) associated with the optimal solution can be interpreted as marginal prices of the commodity at each location in the network. In a similar manner to previous work on auction market formulations for \edit{allocation or scheduling of} compressible natural gas flows \cite{rudkevich2017locational,zlotnik2019optimal}, we propose that the Lagrange multipliers associated with equality constraint \eqref{eq:flowbal} reflect the sensitivity of the objective function to an increment in the withdrawal $d_j$ for $j \in \mathcal N$. Because the objective function $J_P$ in \eqref{eq:combinedobj} that is used in Formulation 3 reflects the total economic value provided to users of the petroleum pipeline system, the Lagrange multipliers that are associated with \eqref{eq:flowbal}, which we denote by $\sigma_j$ for $j\in \mathcal N$, can be interpreted as prices of the transported commodity. In this study we focus on the proposed optimization formulations and the behavior and sensitivity of solutions, which we explore in the section \ref{sec:results}, including an analysis of dual solutions. Further engineering economic analysis of liquid transport pipelines will remain for future work.

\section{Results} \label{sec:results}

We now present a computational study to examine the results of applying the proposed optimization formulations for flow \edit{allocation} and automatic price formation of actual petroleum pipeline systems. The considered test case is synthesized based on physical and economic aspects on the Seaway Pipeline System using openly available information. We estimate unavailable parameters, such as pump properties, and added fictitious supply and consumption points to create a richer variety of possible solutions to examine in this and subsequent studies. We define all of the parameters in the optimization formulation using realistic physical and economic values. 
We perform a sensitivity analysis that focuses on the economic aspects of the pricing bid structure and outcomes. We first briefly review our computational implementation, summarize the solutions for the three formulations, and then discuss the sensitivity analysis outcomes.

\subsection{Computational Implementation}
The optimization problems are solved using IPOPT \cite{wachter2006implementation}, which provides solutions to problems in the standard form \eqref{eq:ipoptform}. The formulations above are implemented using the programming language Julia \cite{Julia1}, an open source and high-performance language for numerical computing. The JuMP package \cite{Julia2} is used as a modeling layer for the optimization problems. JuMP includes a convenient interface to IPOPT and can support a wide range of optimization tasks such as linear programming (LP), nonlinear programming (NLP), and mixed-integer nonlinear programming (MINLP). \edit{The implementation is available  \textsc{PetroleumModels.jl} package \cite{petroleummodels20}, which includes the basis for nonlinear, non-convex problems that optimize a stationary liquid pipeline network, and can be used for academic and research purposes.}

\subsection{Seaway Test Case}
\begin{figure}[t!]
  \centerline{\includegraphics[scale=0.15]{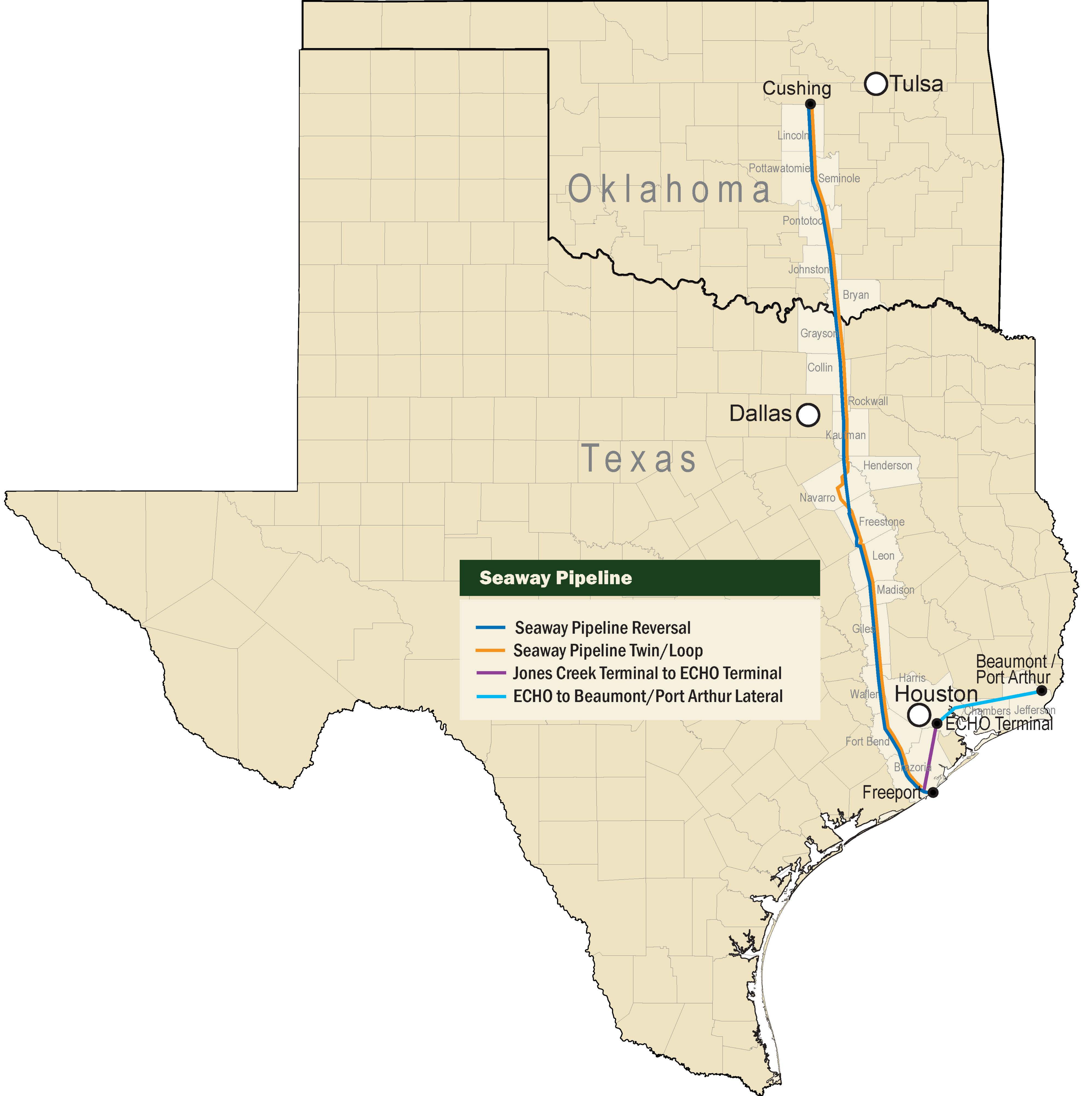}}
  \caption{Seaway Crude Pipeline System (source: https://seawaypipeline.com/)}
  \label{fig:seaway}
  \vspace{-1ex}
  \end{figure}

\begin{figure}[th!]
  \centerline{\includegraphics[width=.75\linewidth]{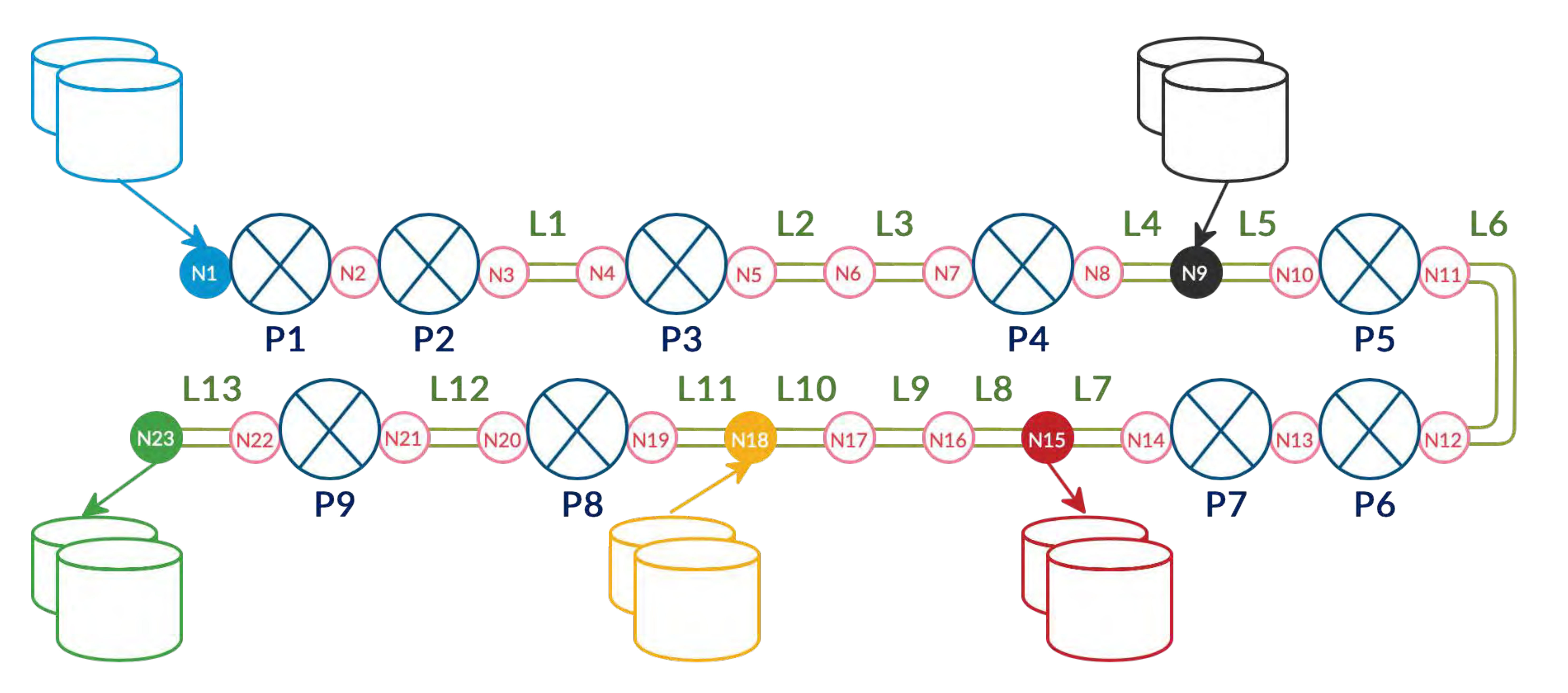}}
  \caption{Schematic of the Seaway Crude Pipeline Test System.  The model has three producers (located at nodes N1, N9, N18) and two consumers (located at nodes N15, N23). There are nine pumps (P1-P9) in the system and thirteen pipes (L1-L13).}
  \label{fig:seawayschmatic}
  \vspace{-1ex}
  \end{figure}

To investigate the performance of algorithms for crude oil transportation problems similar to those proposed in previous studies \cite{13Zhoua, 15Isom, 16Cafaro, 17Wenlong1, 18Wenlong2}, we present a physically and economically realistic pipeline network model and test case that includes auction market parameters drawn from prevailing price and quantity data in January of 2020. The basis for our synthesized example is the Seaway crude oil pipeline system shown in Figure \ref{fig:seaway}, which has an average annual capacity of up to 950,000 barrels per day (bpd), equivalent to 6300 m$^3/$h, for transport through 30 in (76 cm) diameter pipes over a distance of 601 miles (968 km) between Cushing, Oklahoma and Freeport, Texas.
Openly available data includes pipe lengths and diameters, as well as approximate geographical locations, but not metadata regarding hydraulic parameters or pumping equipment. The Seaway system is currently configured to transport crude oil from North to South. A schematic of the Seaway Pipeline System topology is given in Figure \ref{fig:seawayschmatic}. The pipeline delivers crude oil to the Jones Creek terminal, but can also make deliveries to Freeport and Texas docks or to Enterprise ECHO terminal, as well as to Beaumont / Port Arthur, where the Seaway system connects to three terminals: Sunoco Nederland, Phillips 66 Beaumont, and Enterprise Beaumont Marine West.  

  \FloatBarrier

  \FloatBarrier
  We model weakly compressible fluid flow by assuming constant and homogeneous fluid properties for transported crude oil. We use nominal values of 827 kg$/$m$^3$ for density $\rho$, and 4.9$\cdot$10$^{-6}$ m$^2/$s for viscosity $\nu$. The inequality constraints \eqref{eq:pipeflowlim}-\eqref{eq:demandlim} are parameterized by minimum and maximum bound values for pipe, pump, and shipper characteristics as shown in Table \ref{tab:boundval}. The prices offered by the market participants are shown in Table \ref{tab:boundval}, as well. 

\begin{table}[!ht] 
\centering
\renewcommand\thetable{1}
\caption{\textsc{Parameters used in bound values for inequality constraints}}  \label{tab:boundval} 
\begin{threeparttable}
\begin{tabular}{|c|c|c|c|}
\hline
 \multicolumn{4}{|c|}{\textbf{Pipe Inequality Parameters}}\\ \hline \hline
 \thead{} & & \thead{Min value} & \thead{Max value}\\ \hline
 \thead{Flow rate, m$^3$/h} & & {500} & {5000}\\ \hline
  \thead{Head, m} & & {30} & {740}\\ \hline
\multicolumn{4}{|c|}{\textbf{Pump Inequality Parameters}}\\ \hline \hline
 \thead{} &&  \thead{Min value} & \thead{Max value}\\ \hline
 \thead{Efficiency, \%} & &  {60} & {87}\\ \hline
 \thead{Rotation speed, rpm} & & {2400} & {3600}\\ \hline
\multicolumn{4}{|c|}{\textbf{Producer/Consumer Flow Rate Bounds, m$^3/$h}}\\ \hline \hline
 \thead{} & \thead{Bid/Offer \edit{(\$/m$^3$)}} & \thead{Min value} & \thead{Max value}\\ \hline
 \thead{Producer 1 (N1)} & 300 & {360} & {2950}\\ \hline
 \thead{Producer 2 (N9)} & 300 &  {360} & {2950}\\ \hline
  \thead{Producer 3 (N18)} & 300 & {360} & {2950}\\ \hline
  \thead{Consumer 1 (N15)} & 310 & {720} & {3600}\\ \hline
  \thead{Consumer 2 (N23)} & 310 & {720} & {3600}\\ \hline
\end{tabular}
\end{threeparttable}
\end{table}
\FloatBarrier
The lengths of pipeline sections and elevations are shown in Figure \ref{fig:headsol}. In order to complete the hydraulic model, our synthesized network includes pumping machinery with standard equipment parameters (see Table \ref{tab:pumpchar}), which is sufficient to pump the required volume of crude oil at the average statistical rate. 
\begin{table}[h!] 
\centering
\renewcommand\thetable{2}
\caption{\textsc{Standard pump model \& electricity prices used in Seaway test case}}  \label{tab:pumpchar} 
\begin{threeparttable}
\begin{tabular}{|c|c|c|c|c|}
\hline
 \multicolumn{4}{|c|}{\textbf{Nominal parameters}} &\thead{\textbf{H-Q characteristics}}\\ \hline \hline 
 \thead{Flow rate, m$^3/$h} & \thead{Head, m} & \thead{Efficiency, \%} & \thead{Rotational speed, rpm} & $a_0$ = 276.8 \, m\\ \hline
 3600 & 240 & 87 & 3000 & $a_1 = 7.1 \cdot 10^{-6} $\, s$^{2}/$m$^{5}$\\ \hline
\end{tabular}
\begin{tabular}{|c|c|c|c|c|c|c|c|c|c|} \hline
 \thead{Pump ID} & \thead{P1} & \thead{P2} & \thead{P3} & \thead{P4} & \thead{P5} & \thead{P6} &  \thead{P7} & \thead{P8} & \thead{P9} \\ \hline
\thead{\,\,\, Electricity price, \$/kWh \,\,} & {0.12} & {0.12} & {0.13} & {0.11} & {0.12} & {0.15} & {0.15} & {0.08} & {0.14} \\ \hline 
\end{tabular}
\end{threeparttable}
\end{table}
The pumping systems in our model are powered by electric variable frequency drives, for which the variable cost of operation (excluding capital and regular maintenance) is the price of electricity. 
The range of electricity prices was chosen as 0.08-0.15 \$$/$kWh, which corresponds to the prevailing January 2020 price range in the states of Oklahoma and Texas (Table \ref{tab:pumpchar}).
%
%
%
                                                            
\subsection{Optimization Solutions} \label{subsec:solutions}
\begin{figure}[th!]
  \centerline{\includegraphics[width=.6\linewidth]{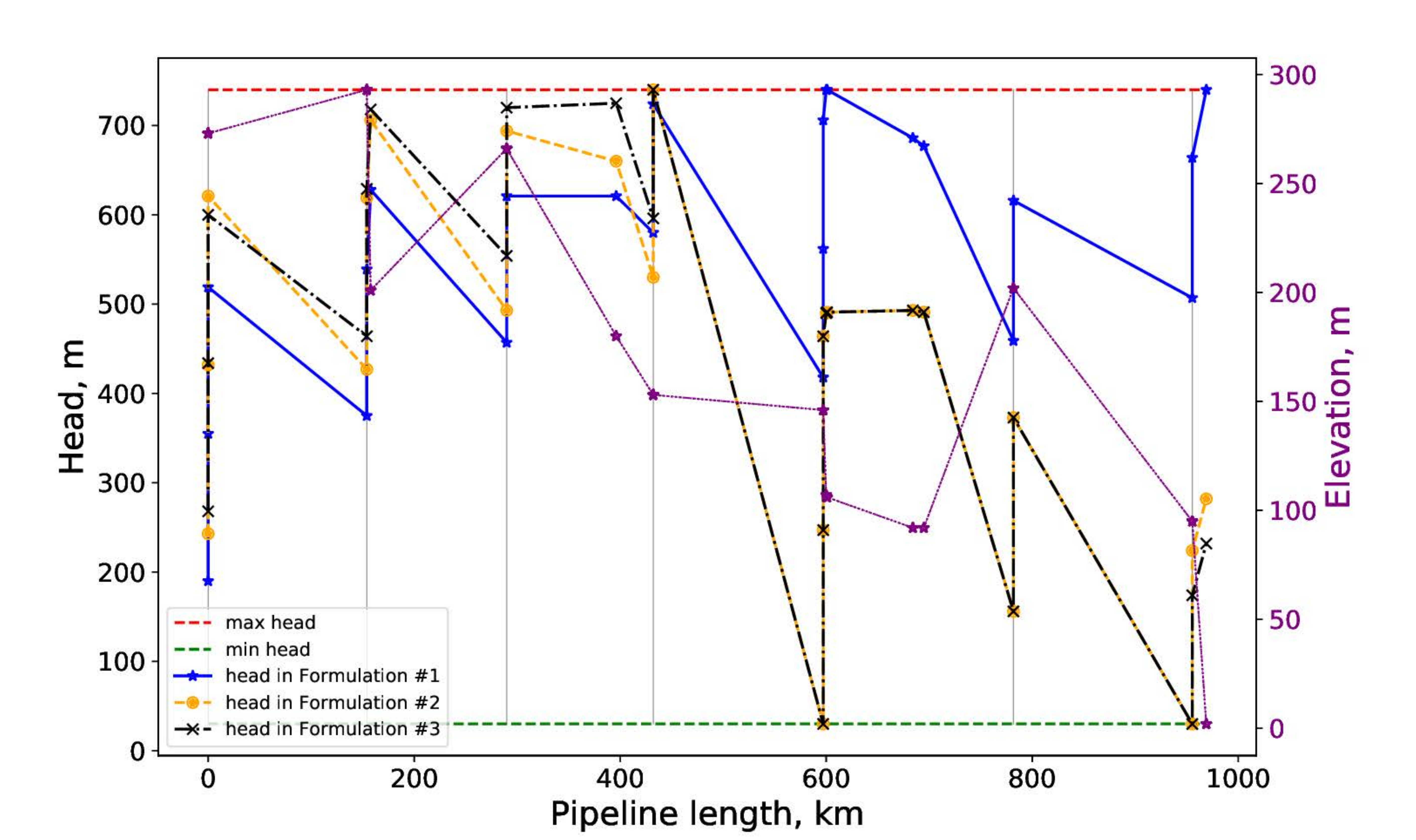}}
  \caption{Elevation along the pipeline (purple) and pump locations (vertical grey lines).  Head solutions are given for Formulations 1 (blue), 2 (yellow), and 3 (black) for the nominal parameter set, as well as lower (green) and upper (red) bounds on head.}
  \label{fig:headsol}
  \vspace{-1ex}
\end{figure}

\begin{figure}[th!]
  \centerline{\includegraphics[width=.6\linewidth]{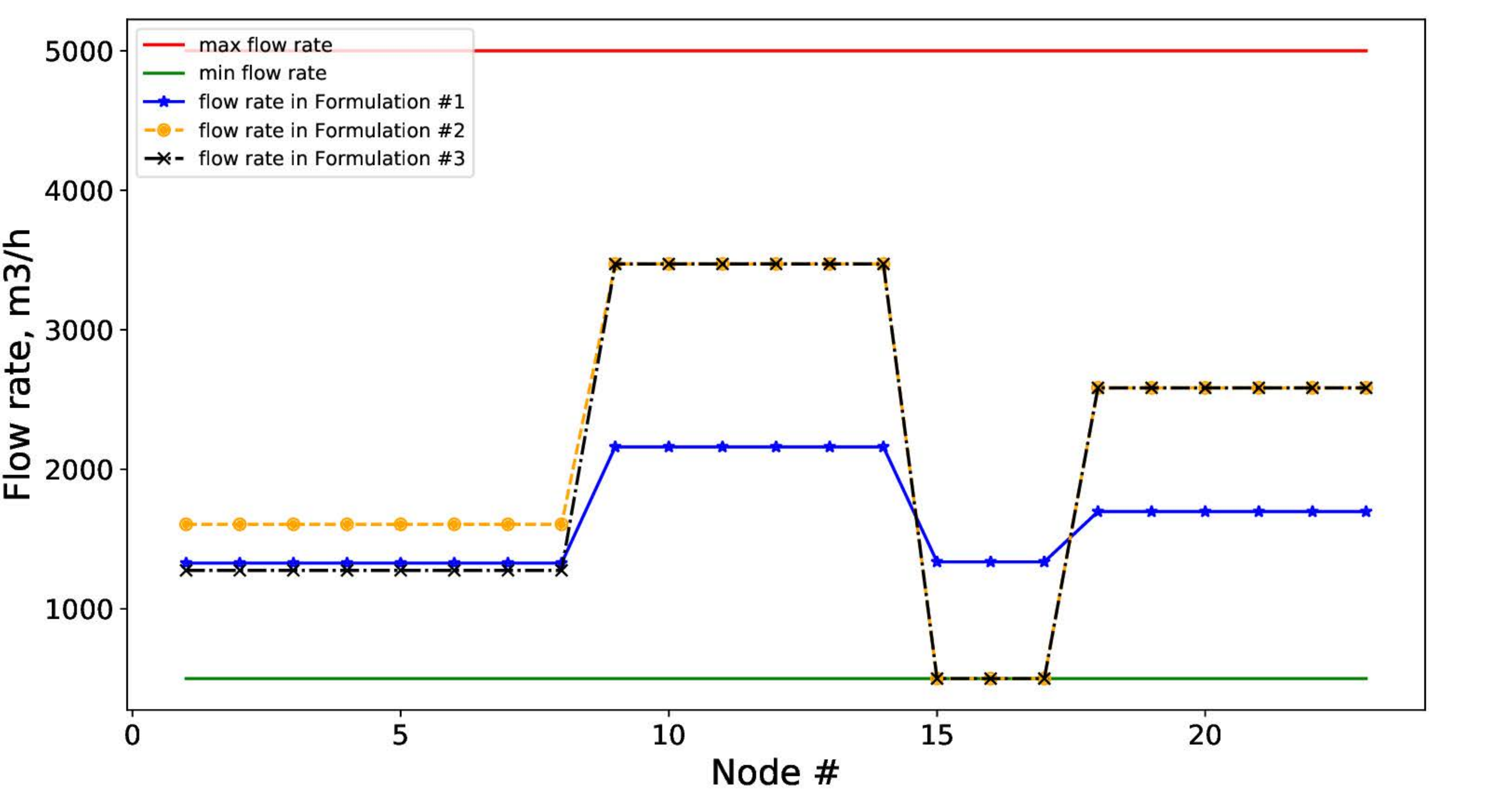}}
  \caption{Flow rates along the pipeline system are shown for the nominal parameter set for solutions to Formulations 1 (blue), 2 (yellow), and 3 (black), and the lower (green) and upper (red) bounds on flow rate are shown.}
  \label{fig:flowsol}
  \vspace{-1ex}
\end{figure}

\begin{figure}[th!]
\centering
\includegraphics[width=.49\linewidth]{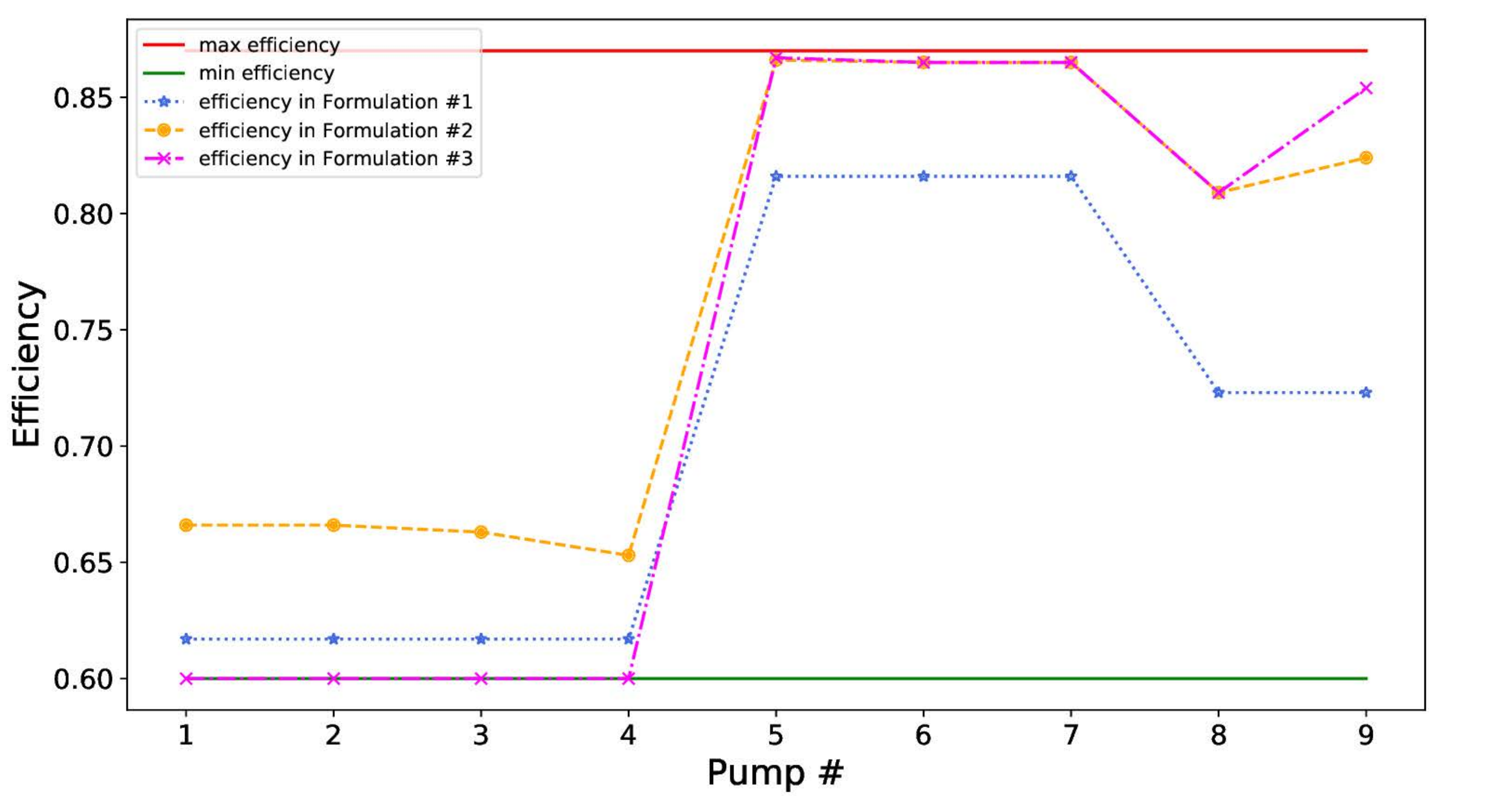}
\includegraphics[width=.49\linewidth]{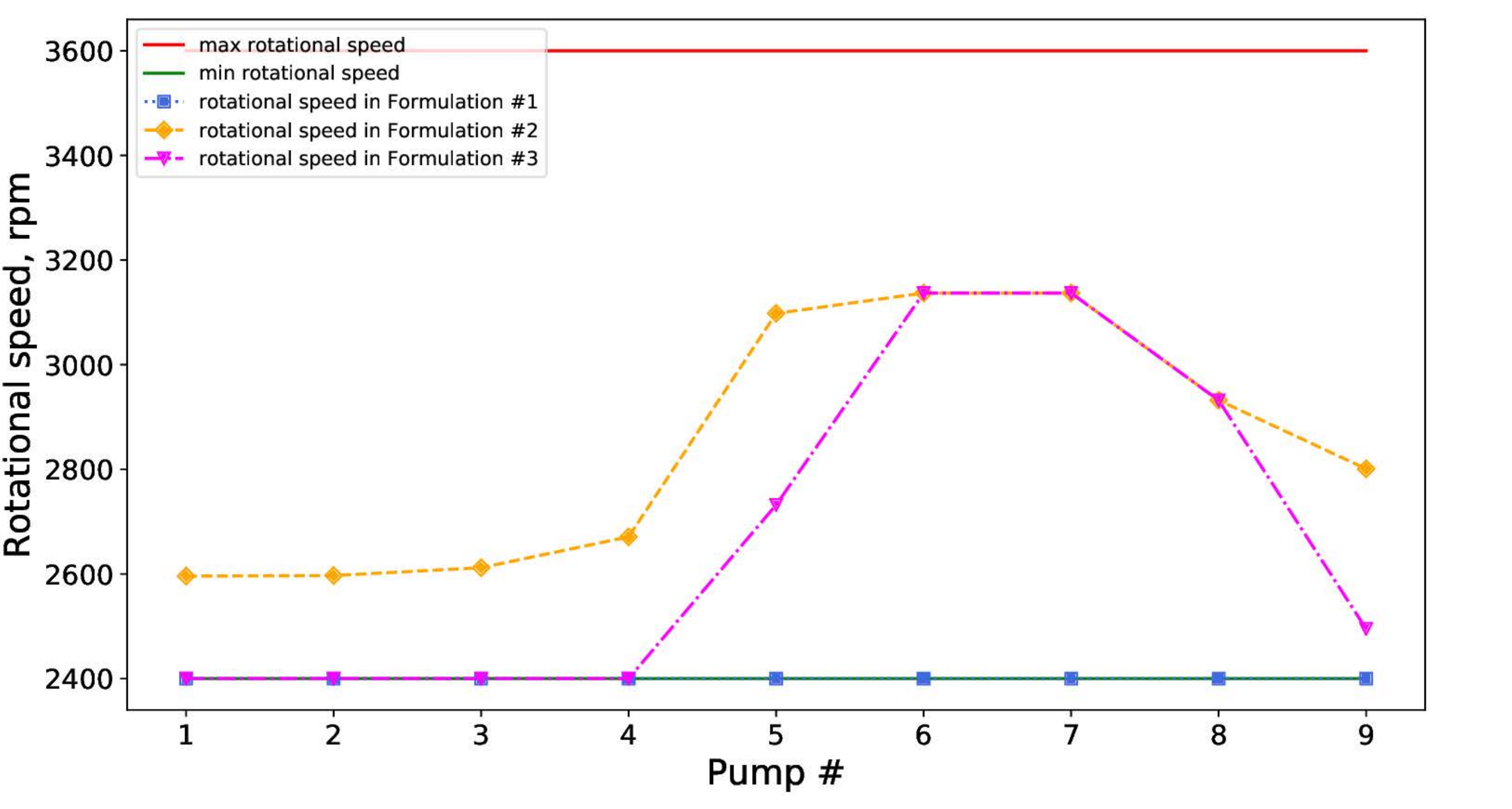}
\caption{Pump efficiency (left) and rotational speed (right) solutions obtained for Formulations 1 (cyan), 2 (yellow), and 3 (magenta) using the baseline parameter set. The lower (green) and upper (red) bounds on both parameters are indicated as well.}
\label{fig:pumpsol}
  \vspace{-1ex}
\end{figure}

The optimization formulations were implemented and solved to local optimality using IPOPT \edit{on a commodity computing platform, which required less than 0.1 seconds in all instances}. The solutions obtained for all three formulations given in Section \ref{subsec:formulations} are presented in Figures \ref{fig:headsol}, \ref{fig:flowsol}, and \ref{fig:pumpsol}. 
\\
\FloatBarrier
First, we use Figure \ref{fig:headsol} to illustrate the geometry (length of pipeline sections) and elevation of the pipeline system, as well as the head configurations in the solutions to the three formulations. The nine main pumps with parameters that are described in Tables \ref{tab:pumpchar} are connected in series, and their locations are indicated with vertical lines in Figure \ref{fig:headsol}.
Observe that because in Formulation 1 the flow rate of consumers is fixed and the optimization objective is to minimize the cost of power used for pump operation, the lower bound constraints on rotation speed are binding for all the pumps.  This results in the least energy system-wide pump operation that is required to satisfy the lower bound constraints on head throughout the system. In contrast, Formulation 2 maximizes the economic value of transportation provided by the system, which results in a redistribution of flows, as seen in Figure \ref{fig:flowsol}.

\noindent Observe that in the solution for Formulation 2, the suppliers at nodes N1 and N9 provide more petroleum, and the first consumer withdraws more than in the solution of Formulation 1. The total volume of the commodity transported in the solution of Formulation 2 is 10.2$\%$ greater than in the solution of Formulation 1. Next, we examine the pump control solutions obtained using the optimization. Figure \ref{fig:pumpsol} shows the efficiency and rotational speed of all nine pumps in the system for each of the three formulations.
Though the increase of 23.35$\%$ in the value of transportation between the solutions to Formulations 1 and 2 comes with an increase of 22.4$\%$ in the cost of operations, the total value created by the system increases by 23.5$\%$. In the solution to Formulation 3, the utilization of capacity on the system is optimized to produce the same flow to consumers and yields the same transportation service value. 
However, the different allocation from the suppliers results in a more efficient pump operating regime. The economic values and costs of transportation and system operations are compared for the three formulations by examining objective function values, as seen in Table \ref{tab:solutions}. 
\begin{table}[h!] 
\centering
\renewcommand\thetable{3}
\caption{\textsc{Objective Function Values for Optimal Solutions}}  
\label{tab:solutions} 
\begin{threeparttable}
\begin{tabular}{|c|c|c|c|}
\hline
\thead{} &\multicolumn{3}{|c|}{\textbf{Expenses values, \$$/$h}}\\ \hline \hline
 \thead{} & {Formulation 1} & {Formulation 2} & {Formulation 3}\\ \hline
  \thead{Transport Value ($J_E$)} & {51484} & \textbf{63509} & {63463}\\ \hline
 \thead{Pumping Energy ($J_O$)} & \textbf{1556} & {1904}& {1742}\\ \hline
  \thead{Total Value ($J_P$)} & {49928} & {61647}& \textbf{61720}\\ \hline
  \thead{Transported Volume ($m^3/h$)}  & {2521} & {5556}& {5556}\\ 
\hline  
\end{tabular}
\end{threeparttable}
\end{table}

In Table \ref{tab:solutions}, the optimized value is highlighted for each formulation. Ultimately, we see that the total economic value created by the pipeline system is maximized by Formulation 3 as well as total transported volume of the product in 2.2 times. Though the cost of pumping energy is 11.9$\%$ greater than in the solution to Formulation 1, the overall increase in economic value is 23.27$\%$. Taking into account an annual pipeline operating time of 8400 hours (with approximately 360 hours set aside for maintenance), the total additional economic value provided by applying the solution to Formulation 3 in comparison with Formulation 1 is $\$100M$  per year, and the total volume of transported commodity increases by over 70\%. Furthermore, if the flow allocation and pumping control cannot be jointly optimized (as in Formulation 3) for business reasons, it is advantageous to prioritize the economic objective function as in Formulation 2 rather than the operating costs as in Formulation 1. Finally, we note that for this baseline parameter set the triples at the bottom of Table \ref{tab:boundval} constitute participant price/quantity bids to sell ($t_j,s_j^{min}, s_j^{max}$) and buy ($r_j,d_j^{min},d_j^{max}$) into the auction market. For these bids, the cleared price solutions $\sigma_j$ at the nodes N1, N9, N18, N15, and N23 with participants are equal to the sell or buy bids $t_j$ or $r_j$, respectively.

\subsection{Sensitivity Analysis} \label{subsec:sensitivity}

Here we discuss a sensitivity analysis of the liquid pipeline flow allocation problem given in Formulation 3, by examining the variation about the nominal solution (obtained with baseline parameters) described in the previous section. Broadly speaking, a sensitivity analysis can assist in determining the set of model parameters that preserve the solution outcome, the types of restrictions on the form of the objective function that yield acceptable results, and the influence of changes to constraint bound values on the objective function value. Further, the ability to change the value of a variable in the optimal solution can be examined by adjusting parameters. In other words, if a binding constraint determines a variable value, determine conditions which cause the constraint to loosen and change the optimal value of the variable.

\begin{figure}[t!]
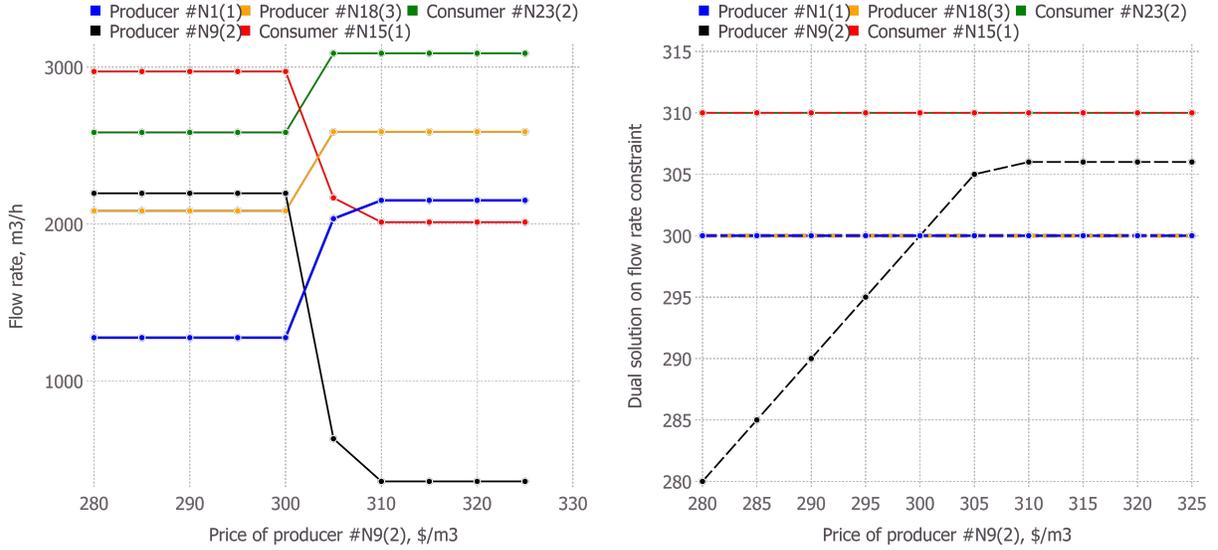

\centering
\includegraphics[width=.49\linewidth]{pic_6.pdf}
\includegraphics[width=.49\linewidth]{pic_3.pdf}
\caption{Sensitivity to the price offered by the second producer of the supply $s_j$ and delivery $d_j$ flow rates (left) and duals $\sigma_j$ of the flow rate constraints (right) in the optimal solution of Formulation 3. Values are shown for producers 1 (node N1, blue), 2 (node N9, black), and 3 (node N18, yellow), and consumers 1 (node N15, red) and 2 (node N23, green).}
\label{fig:sensitivityproducerflowprice}
  \vspace{-1ex}
\end{figure}

The parameters of the pipeline that determine its characteristics both from the economic and technical points of view, are the head developed by the pumping stations, the number of petroleum pumping stations, and pipeline productivity (determined by the flows). All parameters are interconnected, so that a change in one of them entails a change in all the others.  
In this study, we focus on the sensitivity of the physical and marginal pricing outcomes with respect to changes in the price bids $t_j$ and $r_j$ of producers and consumers of the commodity.  
In Figure \ref{fig:sensitivityproducerflowprice}, we show the dependence of supply and delivery flow rates on changes in selling price at which producer 2 (located at node N9). There, we observe that the supply flow rate of producer 2 sharply decreases as the price of producer 2 increases above 300 \$$/$m$^3$. As the price offered by that producer increases past 310, buyers are not willing to purchase the commodity, and the market clears supply by producer 2 at the minimum level required to maintain operating constraints (on head and pumps). The supplies from producers \edit{1 and 3} are increased to compensate, and the demand of consumer 1 decreases. 
Consumer 2 is able to procure more crude oil because producer 3 and consumer 2 are downstream of the price transition at node N9. We also examine the dual solution corresponding to the flow balance constraint, which we interpret as locational marginal prices of crude oil. 
When the offer price of producer 2 is below the baseline of \$300, the price of the commodity at node N9 takes that value. 
However, as the offer price increases past \$300, the cleared market price at node N9 stabilizes at a value of approximately \$306, 
\begin{figure}[!t]
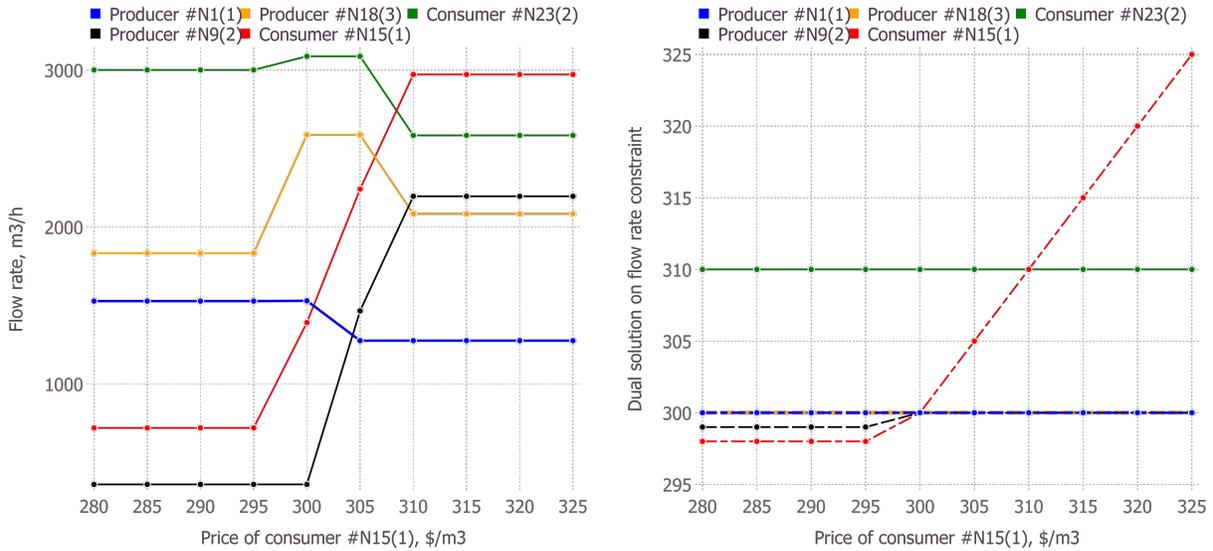

\centering
\includegraphics[width=.49\linewidth]{pic_5.pdf}
\includegraphics[width=.49\linewidth]{pic_4.pdf}
\caption{Sensitivity to the price bid by consumer 1 of the supply $s_j$ and delivery $d_j$ flow rates (left) and duals $\sigma_j$ of the flow rate constraints (right) in the optimal solution of Formulation 3. Values are for producers 1 (node N1, blue), 2 (node N9, black), and 3 (node N18, yellow), and consumers 1 (node N15, red) and 2 (node N23, green).}
\label{fig:sensitivityconsumerflowprice}
\vspace{-1ex}
\end{figure}
which is the maximum value when all other the parameters of the system and participant bids are held at the baseline values.\\ 
Figure \ref{fig:sensitivityconsumerflowprice} is similar to Figure \ref{fig:sensitivityproducerflowprice}, but illustrates the dependence of supply and delivery flow rates and the marginal prices on the bid made by consumer 1. When the bid is below the lowest producer price of \$300, the upstream producers 1 and 2 provide the minimum flow necessary to maintain feasible head and pump characteristics. Note that the cleared value of crude oil for producer 2 and consumer 1 is actually below the bid and offer prices, at \$299 and \$298, respectively. \edit{The cleared price is below this minimum bid because there is a nonlinear component in the objective function that includes the cost of pump operations to bring the commodity downstream. The engineering constraints on the pumps are binding at the Consumer 1 bid values in which this counter-intuitive solution occurs.  The pump efficiency and rotational speed of pumps 1 to 4 bind at the minimum, so that the cleared price reflects the operational requirement to operate these pumps, even though it is not economical to do so.} 
As the bid price of consumer 1 increases above \$300, the incremental value of supply increases for producer 2, so that the flows of both these participants increase significantly until an upper bound constraint on flow in the pipeline binds. 
\begin{figure}[!t]
\centering
\includegraphics[width=.49\linewidth]{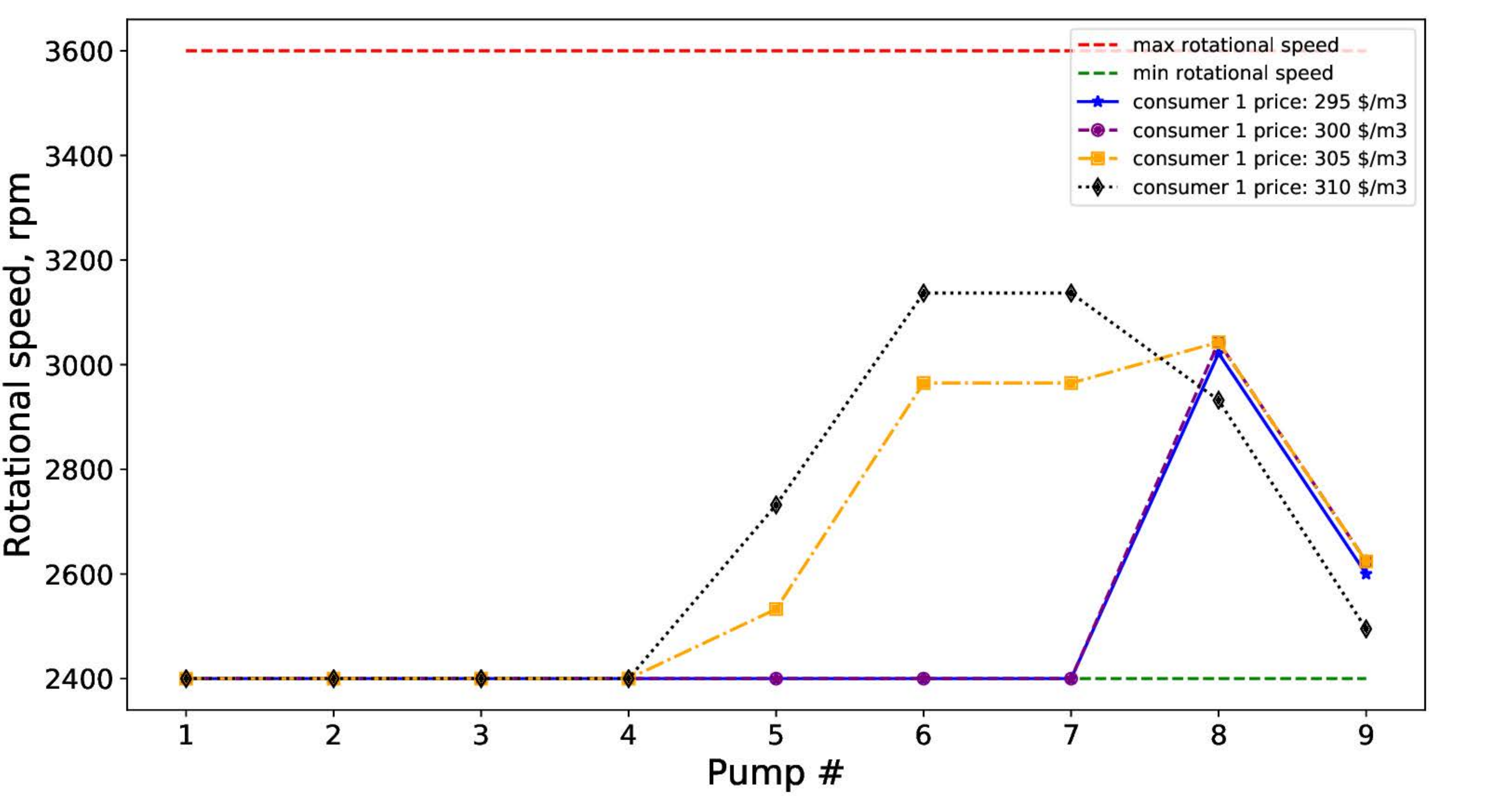}
\includegraphics[width=.49\linewidth]{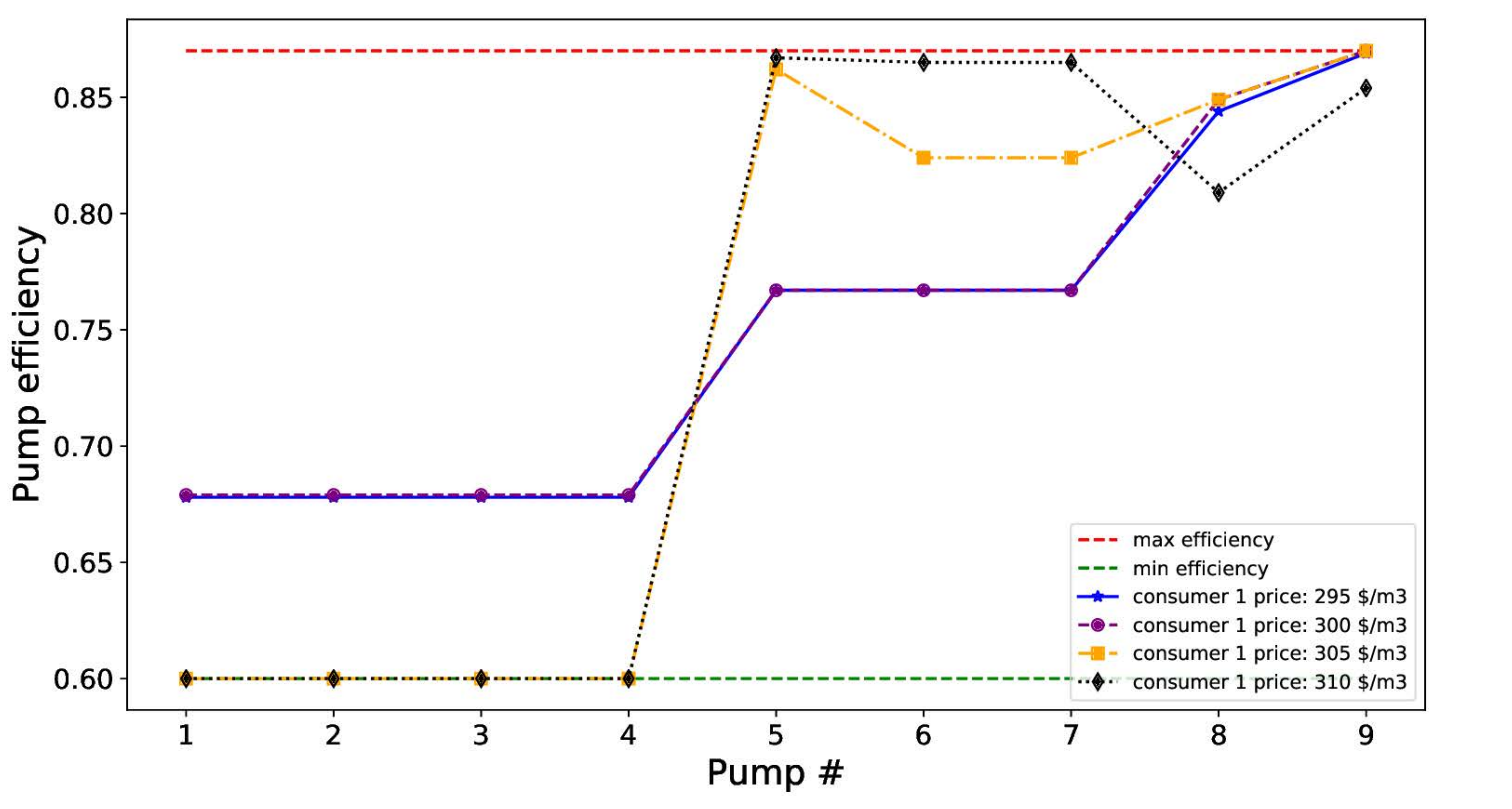} \\
\includegraphics[width=.49\linewidth]{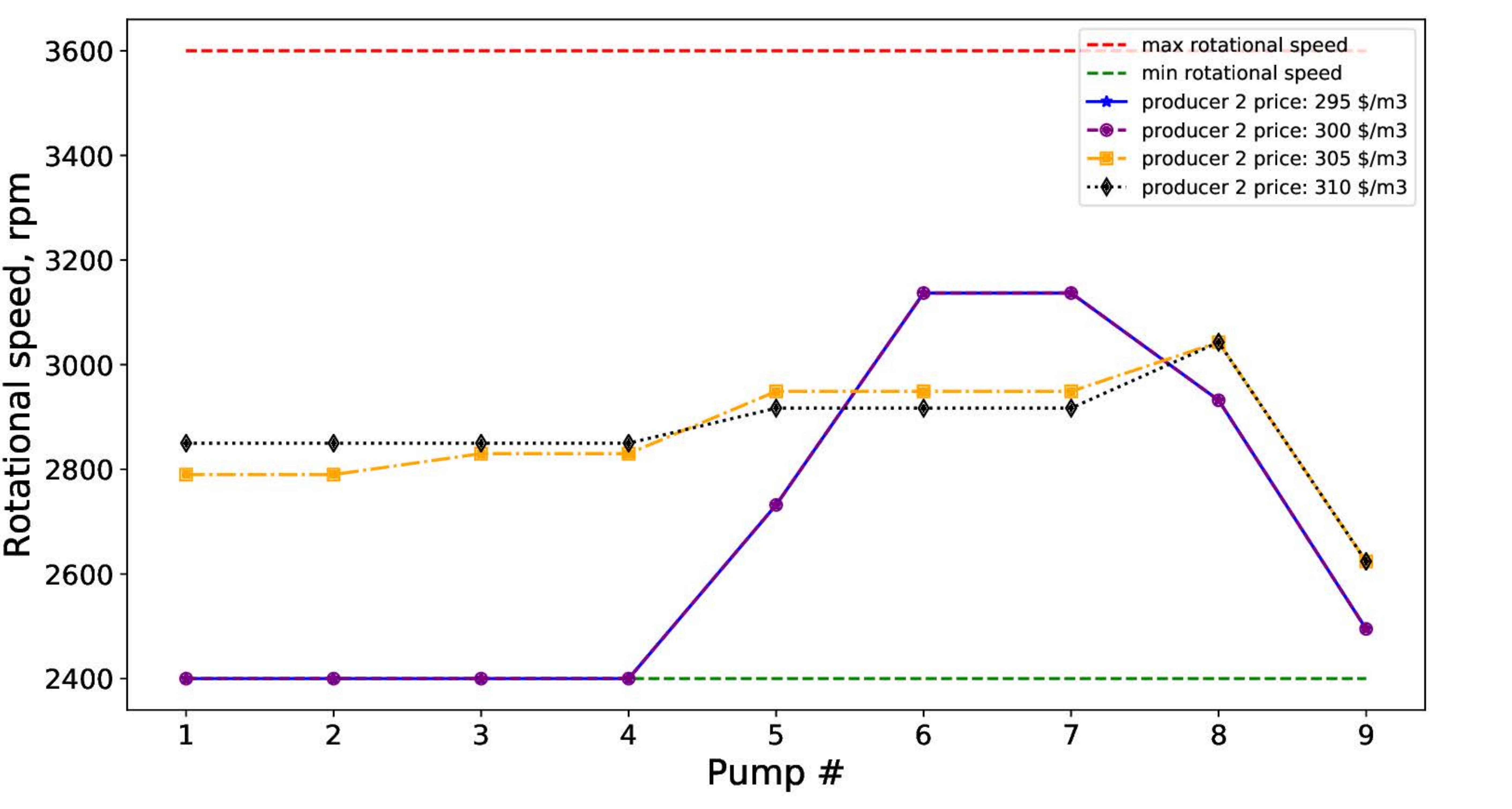}
\includegraphics[width=.49\linewidth]{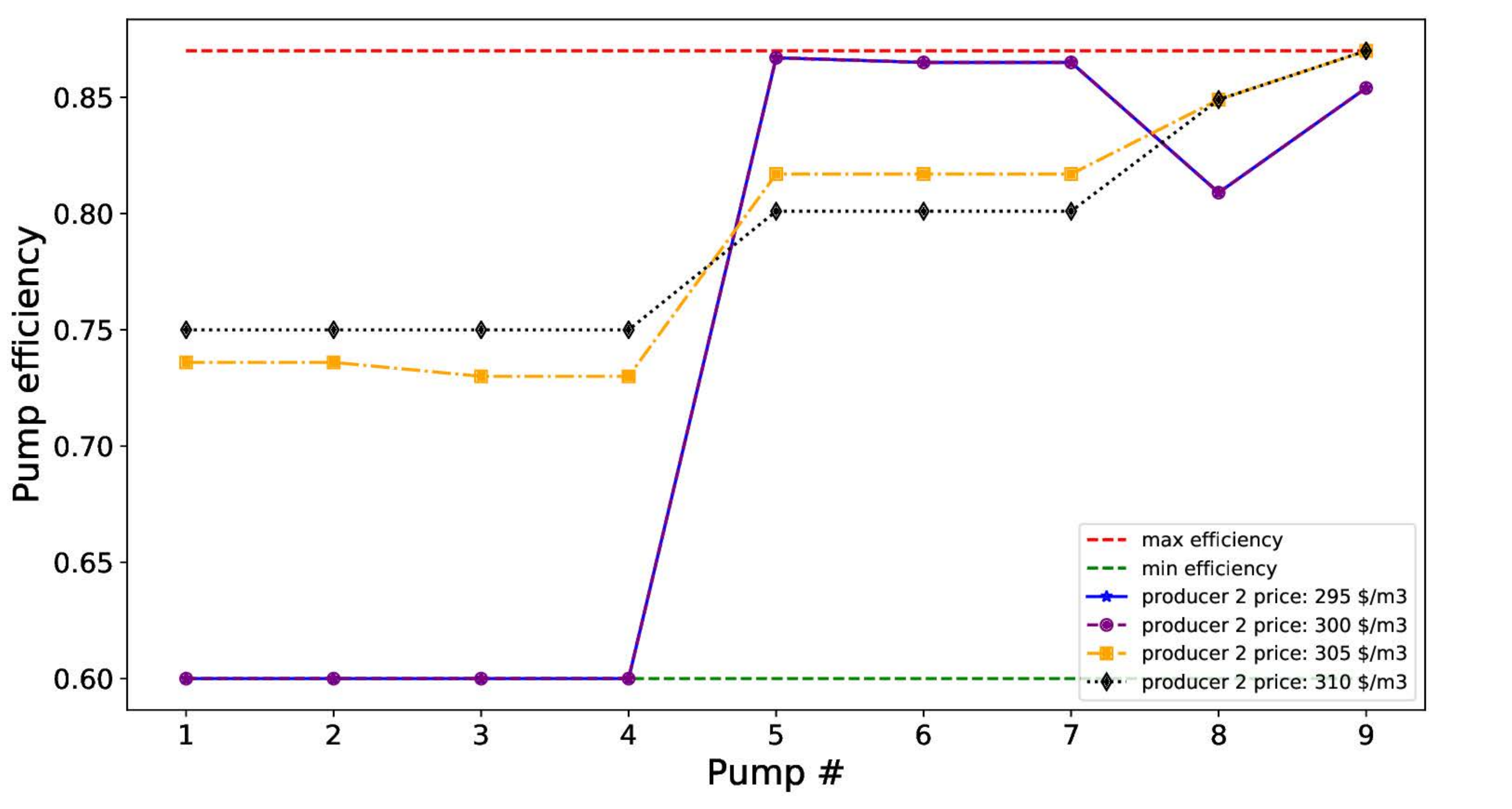}
\caption{Sensitivity to the price bids of the pump rotation speeds and efficiencies in the optimal solution of Formulation 3, with detail near the inflection point. The dependence of rotation speed on the bids of consumer 1 (top left) and producer 2 (bottom left), and the dependence of efficiency on the bids of consumer 1 (top right) and producer 2 (bottom right) are shown. Price levels of \$295 (blue), \$300 (purple), \$305 (yellow), and \$310 (black) are considered, and the upper (red) and lower (green) bound values are indicated.}
\label{fig:sensitivitypumps}
\vspace{-2ex}
\end{figure}
At the consumer 1 price of \$310, the efficiency of pumps P5, P6, and P7 nearly binds at the upper limit (Figure \ref{fig:sensitivitypumps}) resulting in the maximum possible flow of 3472 m$^3/$h on pipe L5. Finally, observe the non-monotone change in the supply flow of producer 3 and delivery flow to consumer 2. This behavior arises from the nonlinear, non-convex nature of the system-wide optimization problem, \edit{and in particular the region defined by the inequality constraint set}. \edit{Such non-monotone properties of sensitivities to optimal solutions in energy networks have been demonstrated in the AC power flow problem \cite{lesieutre2005convexity}.} 
Similar effects can be seen by examining the head, flow, and pump performance configuration of the pipeline system as the price parameters are varied from the inflection point.

\begin{figure}[!ht]
\centering
\includegraphics[width=.49\linewidth]{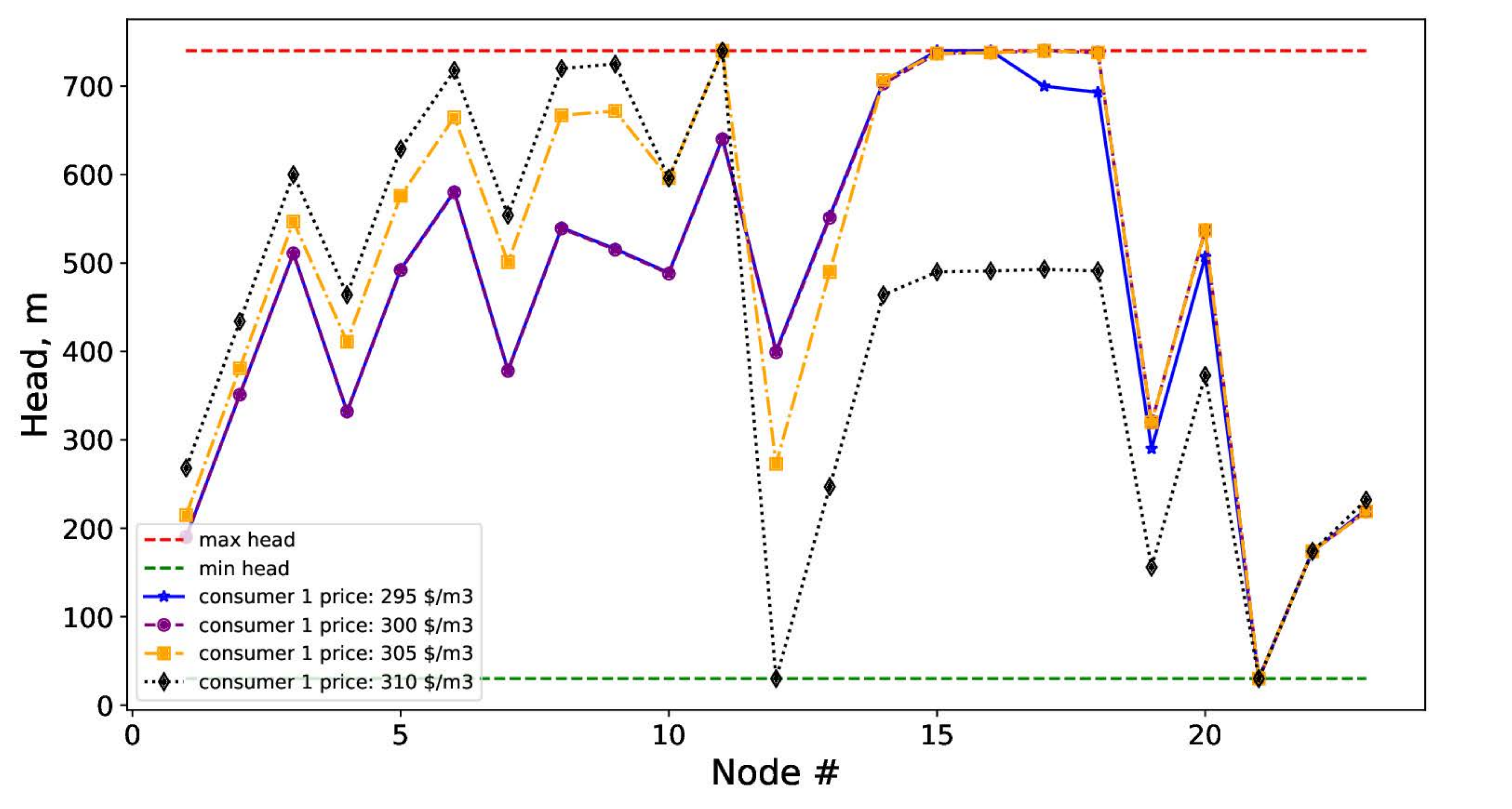}
\includegraphics[width=.49\linewidth]{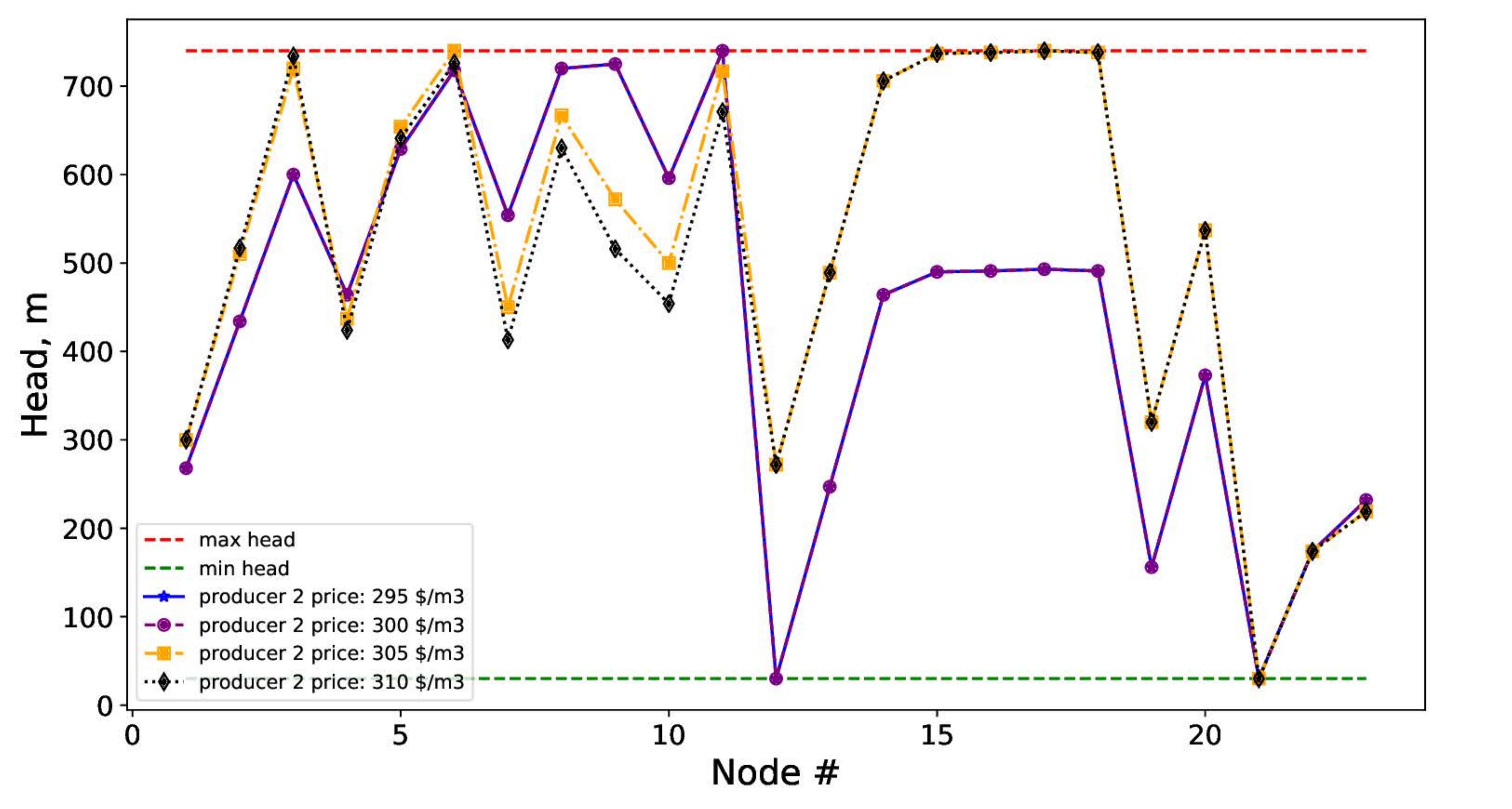}
\caption{Sensitivity to the price bid by the first consumer (left) and the second producer (right) of the head at each node in the optimal solution of Formulation 3, with detail near the inflection point. We consider consumer 1  (producer 2) price bids of \$295 (blue), \$300 (purple), \$305 (yellow), and \$310 (black). The upper (red) and lower (green) bound values on head are shown as well.}
\label{fig:sensitivityconsumerproducerhead}
\vspace{-1ex}
\end{figure}

 \edit{In Figure \ref{fig:sensitivityconsumerproducerhead}, we see that as the price bid of consumer 1 decreases from the nominal value of \$310, consumer 1 is allocated less flow, which causes the head to be decreased upstream of pumps P6 and P7 to maintain head lower bounds, while the head increases downstream because the pumps are able to operate more efficiently. }

Conversely, in Figure \ref{fig:sensitivityconsumerproducerhead}, we see that as the price offer of producer 2 increases from the nominal value of \$300, the flow supplied by producer 2 decreases, so the head upstream of consumer 1 and pumps P6 and P7 decreases, whereas the head downstream of consumer 1 and pumps P6 and P7 increases because the pipeline flow decreases, and the pumps are able to operate more efficiently.  

\begin{figure}[!th]
\centering
\includegraphics[width=.49\linewidth]{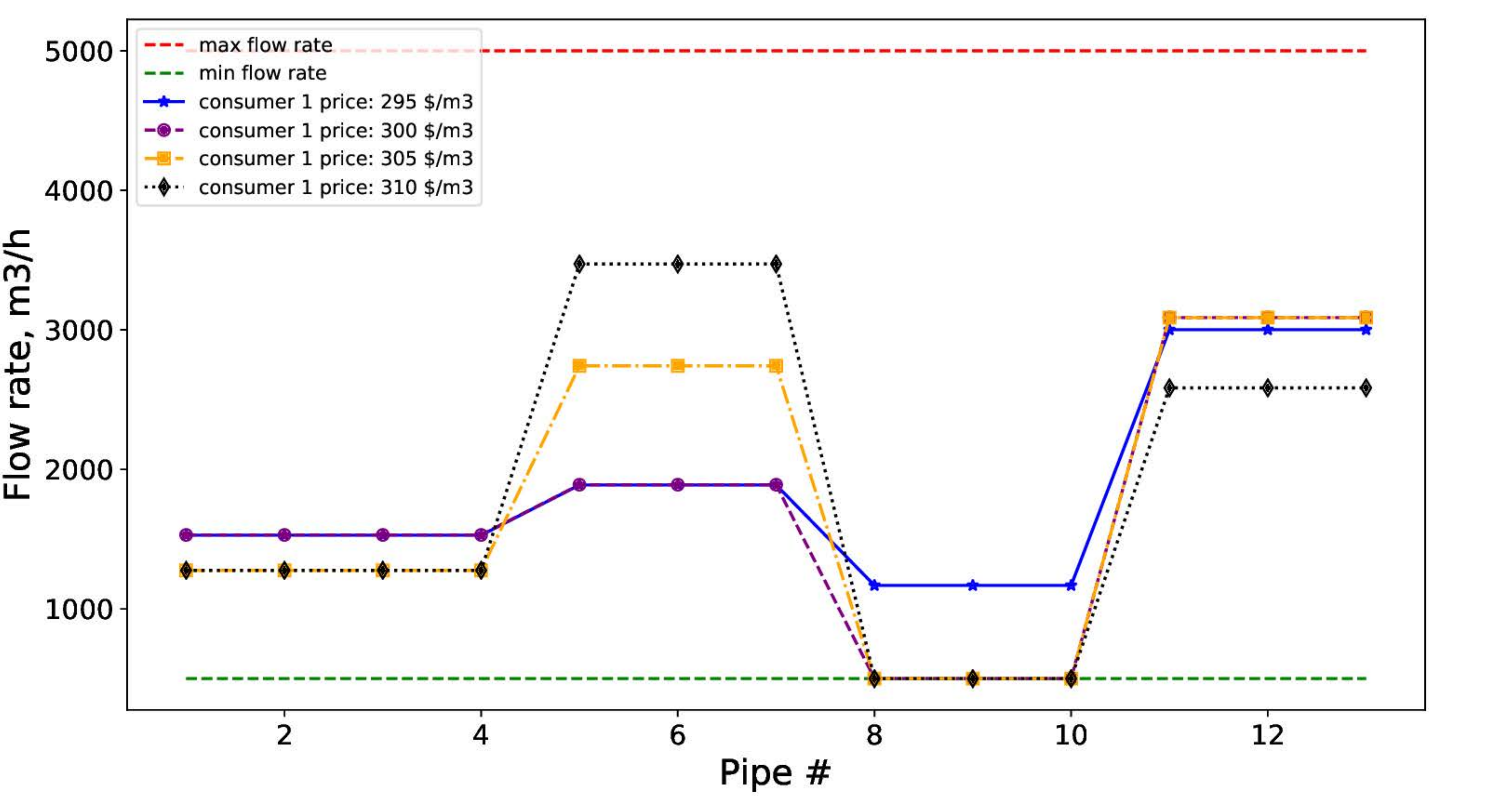}
\includegraphics[width=.49\linewidth]{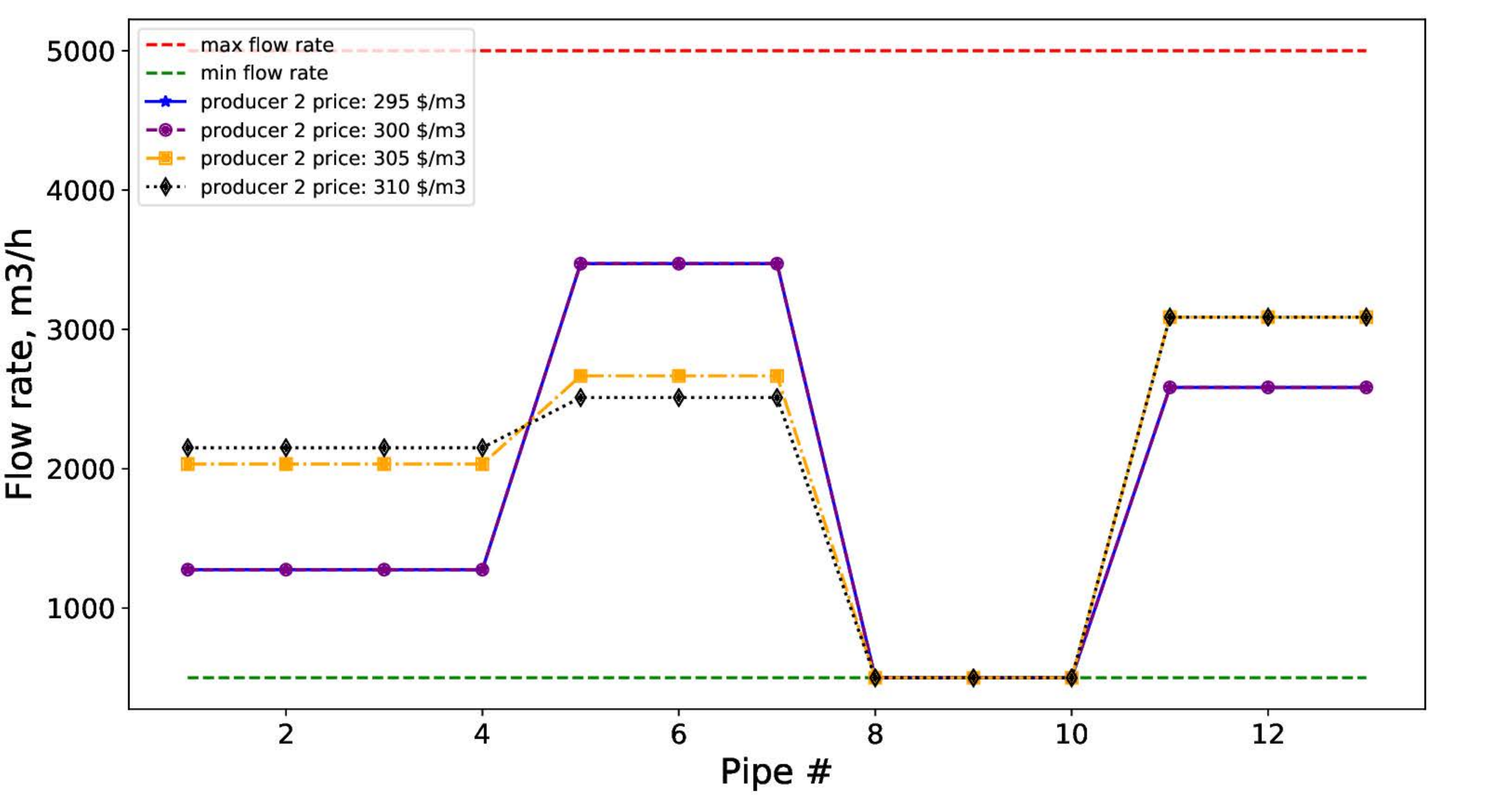}
\caption{Sensitivity to the price bids of the pipeline flow \edit{allocation} in the optimal solution of Formulation 3, with detail near the inflection point. The dependence of flow rates on the bids of consumer 1 (left) and producer 2 (right) at values of \$295 (blue), \$300 (purple), \$305 (yellow), and \$310 (black), as well as the upper (red) and lower (green) bound values on flow, are shown.}
\label{fig:sensitivityflow}
  \vspace{-1ex}
\end{figure}

The relationship between head and flow in the pipeline for optimal solutions under variation of price inputs can be seen in Figure \ref{fig:sensitivityflow}. In the situations where either the consumer 1 price or the producer 2 price is varied, the head in a pipeline section decreases as the flow through that pipeline section decreases, and the converse holds as well. Furthermore, observe in Figure \ref{fig:sensitivitypumps} that pump efficiency increases as rotational speed increases, but also that the pumps operate more efficiently when lower flow rates prevail. The variation in the operating efficiency of the pumps is strongly dependent on location in the pipeline system, and the optimal configuration of head and flow.

\begin{figure}[!th]
\centering
\includegraphics[width=.49\linewidth]{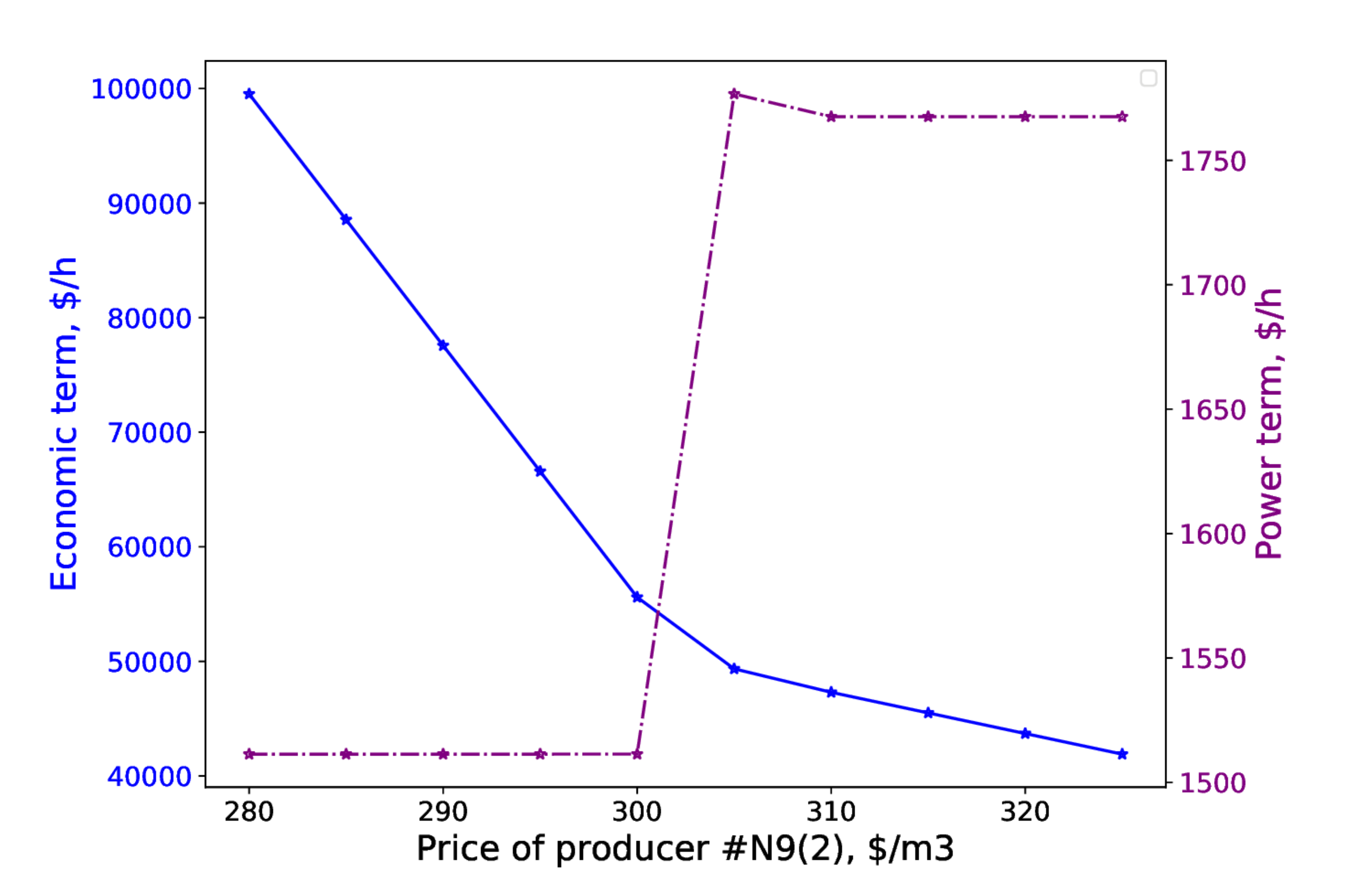}
\includegraphics[width=.49\linewidth]{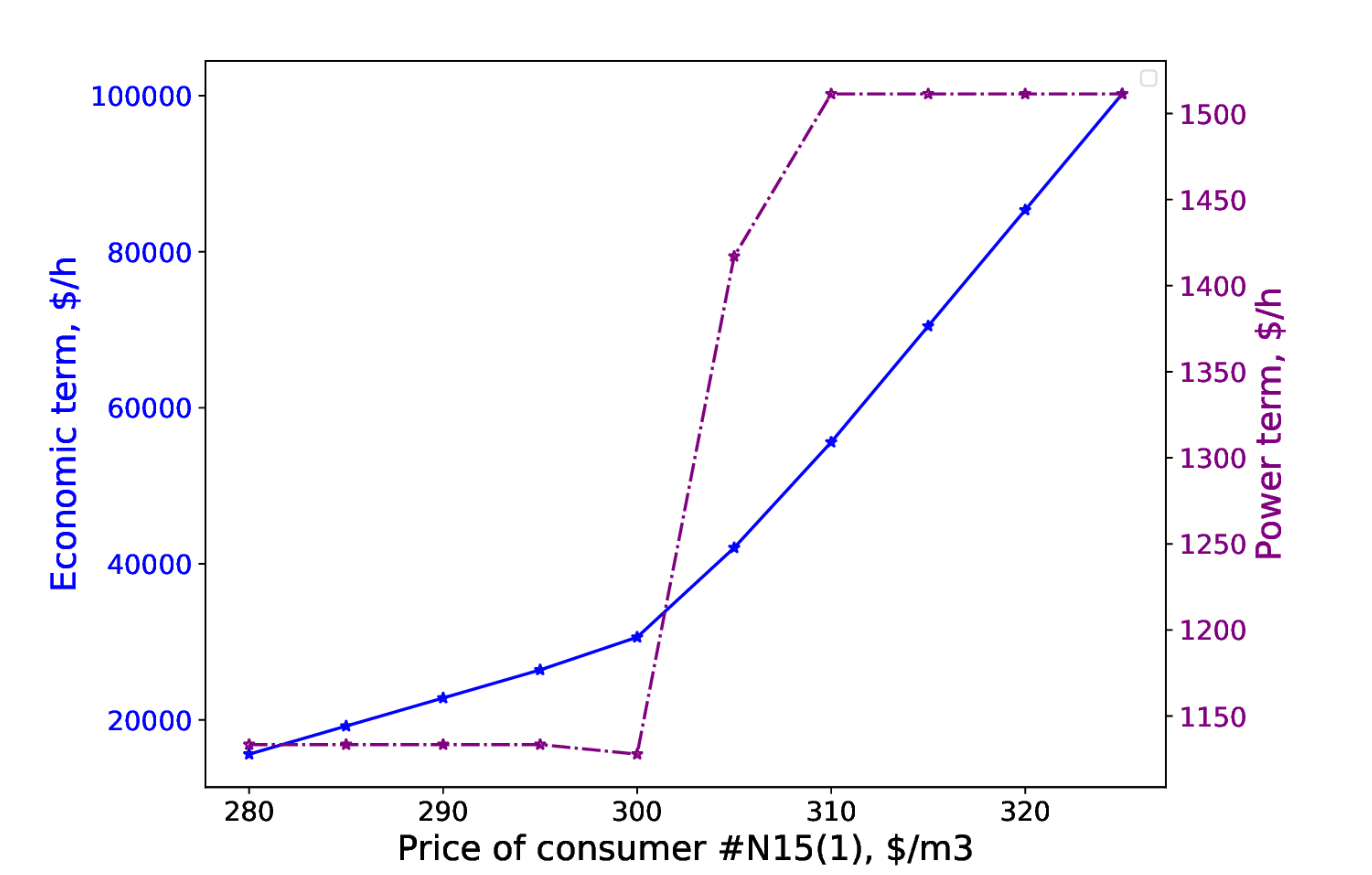}
\caption{Sensitivity to the price bids of the objective function values for the optimal solution of Formulation 3, with detail near the inflection point. The dependencies of the economic component (blue) and the pump operational cost (purple) on producer 2 price (left) and consumer 1 price (right) are shown.}
\label{fig:sensitivityobj}
\end{figure}

Finally, we consider the sensitivity of the terms in the objective function on the variation in the same price parameters as examined above. In Figure \ref{fig:sensitivityobj}, we see that a large economic value is created when the producer price is low, and this value decreases as the producer price increases. A change in the bid by consumer 1 from \$280 to \$325 reflects a significant impact on the economic and energy expenditure components of the objective function, as in the case of the  variation in the offer by producer 2 over the same range of prices. Observe that in both cases, as the consumer bid exceeds the producer offer, the pumping energy is increased to maximize capacity utilization.

\section{Conclusion}

We have developed a general optimization-based methodology for determining flow \edit{allocation} solutions for  liquid pipeline systems \edit{with general network topology} based on engineering economics. The behavior of participants in the petroleum and petroleum products supply chain is modeled using a two-sided single auction market. The problem is solved subject to equality constraints that represent the physics of weakly compressible fluid flow through pipes and models the mechanics of centrifugal pumping equipment that is used to propel fluid flow through the pipeline. The problem includes inequality constraints that reflect the operating limits on flows, pressures, and other state and control variables. Using a realistic case study, we compare and contrast three optimization formulations that minimize operating costs, maximize economic value of transportation, and integrate these first two objectives. Further, we interpret the Lagrange multipliers obtained as optimization output as prices of the transported commodity and perform a perturbation analysis on the economic parameters in order to understand the sensitivity of optimal physical and economic outcomes to changes in the prices that participants bid into the market. We show that substantial economic benefits can be created by the optimal use of frequency-controlled electric drive pumping machinery. With an unregulated electric drive, the efficiency of pumps and electric motors depends on the load and decreases when the pipeline capacity decreases, whereas the degrees of freedom afforded by adjustable electric drives enables enhanced utilization of pipeline capacity. Ultimately, our engineering economics approach could be used by liquid pipeline managers as a decision support tool for integrated financial and physical operations of their systems to be responsive to customer requirements and energy prices. 
The optimization approach that we develop can be calibrated to more specialized situations, such as with specific liquid properties (e.g. viscosity, boundary interactions) \cite{strelnikova2019adaptation},  accounting for the properties of gathering systems \cite{eerc15}, multi-product pipelines \cite{castro2017optimal,relvas2006pipeline,cafaro2009optimal,liang2012hydraulic}, and batching problems \cite{liao2019batch}. Extension to the optimization of water transport systems \cite{hung2011optimization} is possible as well. Fundamentally, the use of more efficient and responsive liquids pipeline \edit{flow allocation and} scheduling could mitigate the substantial carbon dioxide emissions of crude oil transportation and refining \cite{jing2019global}.


\singlespacing
\printbibliography

@book{1lurie,
  title={Modeling of oil product and gas pipeline transportation},
  author={Lurie, M. and Sinaiski, E.},
  pages={234},
  year={2008},
  publisher={Wiley Online Library}
}

@book{2watters,
  title={Modern analysis and control of unsteady flow in pipelines},
  author={Watters, G. },
  year={1980},
  publisher={Ann Arbor Science, Ann Arbor, MI}
}

@book{4thorley2,
  title={Fluid transients in pipeline systems},
  author={Thorley, R.A.D.},
  pages={304},
  year={2004},
  publisher={ASME Press}
}

@book{5chaudhry,
  title={Applied hydraulic transients},
  author={Chaudhry, M.},
  pages={583},
  year={2014},
  publisher={Springer}
}

@book{9korshak2,
  title={Oil, oil products and gas pipeline transportation (in Russian)},
  author={Korshak, A. and Nechval, A.},
  pages={516},
  publisher={Ufa: Design Poligraph Service Publ.},
  year={2005}
}

@article{6lipovka,
  title={Determining hydraulic friction factor for pipeline systems},
  author={Lipovka, A. and Lipovka, Y.},
  journal={Journal of the Siberian Federal University. Engineering \& Technologies},
  pages={62--82},
  year={2014}
}

@article{7Leibenzon_eq,
title={Mathematical modeling of effect of the anti-turbulent additives on flow rate of oil-pipe line section}, 
author={Kuzminskii, Yu. and Shilko, S. and Vyun, V.I.}, 
journal={Journal of Friction and Wear},
  volume={25},
  number={3},
  pages={7--12},
  year={2004}
}

@incollection{brown2003history,
  title={The history of the Darcy-Weisbach equation for pipe flow resistance},
  author={Brown, Glenn O},
  booktitle={Environmental and Water Resources History},
  pages={34--43},
  year={2003}
}

@article{10shabanov1,
  title={Definition of the location of variable frequency drive in the pipeline technological plot (in Russian)},
  author={Shabanov, V. and Bondarenko, O. and Pavlova, Z.},
  journal={Problems of the collection, preparation and transportation of oil petroleum products},
  volume={89},
  number={3},
  pages={93--95},
  publisher={Ufa State Petroleum Technical University, Russia},
  year={2012}
}

@article{11Shabanov2,
   title = {Requirements for the speed of pumps with frequency adjustable electric drive (in Russian)},
    author={Shabanov, V. and Sharipova, S.}, 
    journal={Electrical and data processing facilities and systems},
   volume = {9},
   number = {3},
   year = {2013},
pages={42--46}
}

@article{12Grishin,
title={The efficiency of frequency-regulated electric pump (in Russian)}, 
author={Grishin, V. and Grishin, A.}, 
journal={Automation and Informatization of electrified agricultural production: Scientific papers},
  volume={89},
  pages={118--127},
  publisher={VIESKh, Moscow, Russia},
  year={2004}
}

@article{13Zhoua,
title={Layout optimization of tree-tree gas pipeline network}, 
author={Zhoua, J. and Penga, J. and Lianga, G. and Dengb, T.}, 
journal={Journal of Petroleum Science and Engineering},
  volume={173},
  pages={666--680},
  year={2019}
}

@article{14Samolenkov,
title={Mathematical model of an oil pumping station with variable frequency drive (in Russian)}, 
author={Samolenkov, S. V.}, 
journal={Transport and storage of petroleum products and hydrocarbons},
  number={2},
  pages={21--24},
  year={2013}
}

@article{15Isom,
title={Two methods of data reconciliation for pipeline networks}, 
author={J. D. Isom and A. T. Stamps and A. Esmaili and C. Mancilla}, 
journal={Computers and Chemical Engineering},
  volume={115},
  pages={487--503},
  year={2018}
}

@article{16Cafaro,
title={Rigorous scheduling of mesh-structure refined petroleum pipeline networks}, 
author={Diego C. Cafaro and Jaime Cerdá}, 
journal={Computers and Chemical Engineering},
  volume={38},
  pages={185--203},
  year={2012}
}

@article{17Wenlong1,
title={Adaptive Genetic Algorithm for Steady-State Operation Optimization in Natural Gas Networks}, 
author={Changjun Li and Wenlong Jia and Yi Yang and Xia Wu}, 
journal={Journal of software},
  volume={6},
  number={3},
  pages={452--459},
  year={2011}
}

@article{18Wenlong2,
title={Steady State Modeling and Simulation for Liquid Petroleum Gas Pipeline Networks}, 
author={Wenlong Jia and Changjun Li and Zhuoran Li}, 
journal={Proceedings of the Pipelines 2013 Conference},
pages={1502--1511},
year={2013}
}

@article {barkhatov2017development,
   title = {Development of methods for energy-efficient operation of trunk pipelines based on the optimization of technological regimes},
   author = {Barkhatov, A. F.},
   journal = {Gubkin Russian State University of Oil and Gas},
   year = {2017}
}

@inproceedings{wang2012survey,
  title={A survey on oil/gas pipeline optimization: Problems, methods and challenges},
  author={Wang, Yu and Tian, Chun-Hua and Yan, Jun-chi and Huang, Jin},
  booktitle={Proceedings of 2012 IEEE International Conference on Service Operations and Logistics, and Informatics},
  pages={150--155},
  year={2012},
  organization={IEEE}
}

@article{relvas2006pipeline,
  title={Pipeline scheduling and inventory management of a multiproduct distribution oil system},
  author={Relvas, Susana and Matos, Henrique A and Barbosa-P{\'o}voa, Ana Paula FD and Fialho, Joao and Pinheiro, Ant{\'o}nio S},
  journal={Industrial \& Engineering Chemistry Research},
  volume={45},
  number={23},
  pages={7841--7855},
  year={2006},
  publisher={ACS Publications}
}

@misc{petroleummodels20,
  author={Khlebnikova, Elena and Zlotnik, Anatoly and Sundar, Kaarthik and Ewers, Mary and Tasseff, Byron and Bent, Russell},
  note={\url{https://github.com/lanl-ansi/PetroleumModels.jl/}},
  year={2020}
}

@inproceedings{khlebnikova2020optimization,
  title={Optimization of Liquid Pipeline Control for Economic and Efficient Operations},
  author={Khlebnikova, Elena and Zlotnik, Anatoly and Sundar, Kaarthik and Ewers, Mary and Tasseff, Byron and Bent, Russell},
  booktitle={SPE Europec featured at 82nd EAGE Conference and Exhibition},
  year={2020},
  organization={Society of Petroleum Engineers}
}

@article{hung2011optimization,
  title={Optimization of water systems with the consideration of pressure drop and pumping},
  author={Hung, Szu Wen and Kim, Jin-Kuk},
  journal={Industrial \& Engineering Chemistry Research},
  volume={51},
  number={2},
  pages={848--859},
  year={2011},
  publisher={ACS Publications}
}

@article{reddy2004novel,
  title={Novel solution approach for optimizing crude oil operations},
  author={Reddy, P Chandra Prakash and Karimi, IA and Srinivasan, R},
  journal={AIChE Journal},
  volume={50},
  number={6},
  pages={1177--1197},
  year={2004},
  publisher={Wiley Online Library}
}

@article{zhou2015dynamic,
  title={Dynamic optimization of heated oil pipeline operation using PSO--DE algorithm},
  author={Zhou, Ming and Zhang, Yu and Jin, Shijiu},
  journal={Measurement},
  volume={59},
  pages={344--351},
  year={2015},
  publisher={Elsevier}
}

@article{liu2014optimal,
  title={Optimal energy consumption analysis of natural gas pipeline},
  author={Liu, Enbin and Li, Changjun and Yang, Yi},
  journal={The Scientific World Journal},
  volume={2014},
  year={2014},
  publisher={Hindawi}
}

@article{liu2015research,
  title={Research on the optimal energy consumption of oil pipeline},
  author={Liu, Enbin and Li, Changjun and Yang, Liuting and Liu, Song and Wu, Mingchang and Wang, Di},
  journal={Journal of environmental biology},
  volume={36},
  number={4},
  pages={703},
  year={2015}
}

@book{wylie1993fluid,
  title={Fluid transients in systems},
  author={Wylie, E Benjamin and Streeter, Victor Lyle and Suo, Lisheng},
  volume={1},
  year={1993},
  publisher={Prentice Hall Englewood Cliffs, NJ}
}

@article{losenkov2019optimization,
  title={Optimization of Oil-Flow Scheduling in Branched Pipeline Systems},
  author={Losenkov, A. S. and Yushchenko, T. S. and Strelnikova, S. A. and Michkova, D.E.},
  journal={Journal of Pipeline Systems Engineering and Practice},
  volume={10},
  number={3},
  pages={04019014},
  year={2019},
  publisher={American Society of Civil Engineers}
}

@inproceedings{strelnikova2019adaptation,
  title={Adaptation of fluid motion mathematical model in pipelines using drag reducing agents},
  author={Strelnikova, S. and Yushchenko, T.},
  booktitle={PSIG Annual Meeting},
  year={2019},
  organization={Pipeline Simulation Interest Group}
}

@inproceedings{losenkov2019mathematical,
  title={A mathematical approach to optimization of oil flow scheduling in pipeline systems},
  author={Losenkov, A. and Yushchenko, T. and Strelnikova, S. and Michkova, D.},
  booktitle={PSIG Annual Meeting},
  year={2019},
  organization={Pipeline Simulation Interest Group}
}

@article{mirhassani2008operational,
  title={An operational planning model for petroleum products logistics under uncertainty},
  author={MirHassani, SA},
  journal={Applied Mathematics and Computation},
  volume={196},
  number={2},
  pages={744--751},
  year={2008},
  publisher={Elsevier}
}

@inproceedings{barreto2004optimization,
  title={Optimization of pump energy consumption in oil pipelines},
  author={Barreto, Claudio Veloso and Pires, Luis Fernando Gonc and Azevedo, Luis Fernando Alzuguir and others},
  booktitle={2004 International Pipeline Conference},
  pages={23--27},
  year={2004},
  organization={American Society of Mechanical Engineers Digital Collection}
}

@article{liang2012hydraulic,
  title={Hydraulic model optimization of a multi-product pipeline},
  author={Liang, Yongtu and Li, Ming and Li, Jiangfei},
  journal={Petroleum Science},
  volume={9},
  number={4},
  pages={521--526},
  year={2012},
  publisher={Springer}
}

@article{cafaro2011detailed,
  title={Detailed scheduling of operations in single-source refined products pipelines},
  author={Cafaro, Vanina G and Cafaro, Diego C and M{\'e}ndez, Carlos A and Cerd{\'a}, Jaime},
  journal={Industrial \& Engineering Chemistry Research},
  volume={50},
  number={10},
  pages={6240--6259},
  year={2011},
  publisher={ACS Publications}
}

@article{cafaro2015minlp,
  title={MINLP model for the detailed scheduling of refined products pipelines with flow rate dependent pumping costs},
  author={Cafaro, VG and Cafaro, DC and M{\'e}ndez, CA and Cerd{\'a}, J},
  journal={Computers \& Chemical Engineering},
  volume={72},
  pages={210--221},
  year={2015},
  publisher={Elsevier}
}

@article{mostafaei2017continuous,
  title={Continuous-time scheduling formulation for straight pipelines},
  author={Mostafaei, H and Castro, PM},
  journal={AIChE Journal},
  volume={63},
  number={6},
  pages={1923--1936},
  year={2017},
  publisher={Wiley Online Library}
}

@inproceedings{magalhaes2003crude,
  title={Crude oil scheduling},
  author={Magalhaes, Marcus V and Shah, Nilay},
  booktitle={Proceedings of the 4th Conference on Foundations of Computer-Aided Process Operations},
  pages={323--326},
  year={2003}
}

@article{cafaro2015optimization,
  title={Optimization model for the detailed scheduling of multi-source pipelines},
  author={Cafaro, Vanina G and Cafaro, Diego C and M{\'e}ndez, Carlos A and Cerd{\'a}, Jaime},
  journal={Computers \& Industrial Engineering},
  volume={88},
  pages={395--409},
  year={2015},
  publisher={Elsevier}
}

@inproceedings{starikov2015minimization,
  title={Minimization of pump energy losses in dynamic automatic control of pressure in the main oil pipeline},
  author={Starikov, DP and Rybakov, EI and Gromakov, Evgeniy Ivanovich and Zamyatina, OM},
  booktitle={Future Communication, Information and Computer Science: Proceedings of the 2014 International Conference on Future Communication, Information and Computer Science (FCICS 2014), May 22-23, 2014, Beijing, China.},
  pages={69},
  year={2015},
  organization={CRC Press}
}

@inproceedings{rizwan2013crude,
  title={Crude Oil Network Modeling, Simulation and Optimization: Novel Approach and Operational Benefits},
  author={Rizwan, Mohamed and Al-Otaibi, Mohammed F and Al-Khaledi, Sadoun},
  booktitle={ASME 2013 India Oil and Gas Pipeline Conference},
  pages={V001T02A007--V001T02A007},
  year={2013},
  organization={American Society of Mechanical Engineers}
}

@article{soud14agent,
  title={Agent-Based Simulation Model for Petroleum Supply Chain Management (In Russian)},
  author={Soud, Abdalases Mohammed Amen},
  journal={Information Processing Systems},
  number={9},
  pages={199--205},
  year={2014},
  publisher={V. N. Karazin Kharkiv National University (Ukraine)}
}

@article{camacho1990optimal,
  title={Optimal operation of pipeline transportation systems},
  author={Camacho, EF and Ridao, MA and Ternero, JA and Rodriguez, JM},
  journal={IFAC Proceedings Volumes},
  volume={23},
  number={8},
  pages={455--460},
  year={1990},
  publisher={Elsevier}
}

@book{zemenkov2017,
   title = {Oil and Gas Pipeline and Product Pipeline Operations Engineer Handbook},
   author = {Zemenkov, Y.},
   year = {2017},
   publisher = {Library of Oil and Gas Producers and Contractors}
}

@book{eerc15,
   title = {Liquids Gathering Pipelines: A Comprehensive Analysis},
   author = {Energy \& Environmental Research Center},
   year = {2017},
   publisher = {University of North Dakota}
}

@article{castro2017optimal,
  title={Optimal scheduling of multiproduct pipelines in networks with reversible flow},
  author={Castro, Pedro M},
  journal={Industrial \& Engineering Chemistry Research},
  volume={56},
  number={34},
  pages={9638--9656},
  year={2017},
  publisher={ACS Publications}
}

@article{cafaro2009optimal,
  title={Optimal scheduling of refined products pipelines with multiple sources},
  author={Cafaro, Diego C and Cerd{\'a}, Jaime},
  journal={Industrial \& Engineering Chemistry Research},
  volume={48},
  number={14},
  pages={6675--6689},
  year={2009},
  publisher={ACS Publications}
}

@techreport{pharris2008argonne,
  title={Overview of the design, construction, and operation of interstate liquid petroleum pipelines},
  author={Pharris, TC and Kolpa, RL},
  year={2008},
  institution={Argonne National Lab.(ANL), Argonne, IL (United States)}
}

@article {losenokov2017optimization,
   title = {Optimization of planning of oil cargo flows in branched systems of oil trunk pipelines (in Russian)},
   author = {Losenkov, Alexander Stanislavovich and Yushchenko, Taras Sergeevich and Strelnikova, Svetlana Alekseevna and Michkova, Diana Evgenievna},
   journal = {Industrial Automation},
   number = {8},
   pages = {16--22},
   year = {2017}
}

@article {Julia1,
   title = {Julia: A Fast Dynamic Language for Technical Computing},
   author = {Bezanson, J. and Karpinski, S. and Shah, V. B. and Edelman, A.},
   journal = {[Online]. Available: https://arxiv.org/abs/1209.5145},
   volume = {abs/1209.5145},
   year = {2012}
}

@article {Julia2,
   title = {“JuMP: A modeling language for mathematical optimization},
   author = {Dunning, I. and Huchette, J. and Lubin, M.},
   journal = {SIAM Review},
   volume = {59},
   number = {2},
   pages = {295--320},
   year = {2017}
}

@article{wachter2006implementation,
  title={On the implementation of an interior-point filter line-search algorithm for large-scale nonlinear programming},
  author={W{\"a}chter, Andreas and Biegler, Lorenz T},
  journal={Mathematical programming},
  volume={106},
  number={1},
  pages={25--57},
  year={2006},
  publisher={Springer}
}

@article{moro2004mixed,
  title={Mixed-integer programming approach for short-term crude oil scheduling},
  author={Moro, Lincoln FL and Pinto, Jos{\'e} M},
  journal={Industrial \& engineering chemistry research},
  volume={43},
  number={1},
  pages={85--94},
  year={2004},
  publisher={ACS Publications}
}

@article{rejowski2008novel,
  title={A novel continuous time representation for the scheduling of pipeline systems with pumping yield rate constraints},
  author={Rejowski Jr, R and Pinto, Jos{\'e} M},
  journal={Computers \& Chemical Engineering},
  volume={32},
  number={4-5},
  pages={1042--1066},
  year={2008},
  publisher={Elsevier}
}

@article{sarbu1998energetic,
  title={Energetic optimization of water pumping in distribution systems},
  author={S{\'a}rbu, Ioan and Borza, Ioan},
  journal={Periodica Polytechnica Mechanical Engineering},
  volume={42},
  number={2},
  pages={141--152},
  year={1998}
}

@article{georgescu2014estimation,
  title={Estimation of the efficiency for variable speed pumps in EPANET compared with experimental data},
  author={Georgescu, A-M and Cosoiu, C-I and Perju, S and Georgescu, S-C and Hasegan, L and Anton, A},
  journal={Procedia Engineering},
  volume={89},
  pages={1404--1411},
  year={2014},
  publisher={Elsevier}
}

@article{kocher2016analytical,
  title={Analytical Solution for Speed to Achieve a Desired Operating Point for a Given Fan or Pump},
  author={Kocher, Michael F and Subbiah, Jeyamkondan},
  journal={Applied Engineering in Agriculture},
  volume={32},
  number={6},
  pages={751--757},
  year={2016},
  publisher={American Society of Agricultural and Biological Engineers}
}

@article{2012Vyazunov,
  title={Analytical representation of characteristics of centrifugal pumps under variable rotation speed (in Russian)},
  author={Vyazunov E.V., Putin S.V.},
  journal={Science and technologies: oil and oil products pipeline transportation},
  number={4},
  pages={67--69},
  year={2012}
}

@article{marchi2012assessing,
  title={Assessing variable speed pump efficiency in water distribution systems},
  author={Marchi, Angela and Simpson, Angus R and Ertugrul, Nesimi},
  journal={Drinking Water Engineering and Science},
  volume={5},
  number={1},
  pages={15--21},
  year={2012},
  publisher={Copernicus GmbH}
}

@article{strogen2016environmental,
  title={Environmental, public health, and safety assessment of fuel pipelines and other freight transportation modes},
  author={Strogen, Bret and Bell, Kendon and Breunig, Hanna and Zilberman, David},
  journal={Applied Energy},
  volume={171},
  pages={266--276},
  year={2016},
  publisher={Elsevier}
}

@article{jing2019global,
  title={Global crude oil refining: carbon intensity and mitigation potential},
  author={Jing, Liang and El-Houjeiri, Hassan M and Monfort, Jean-Christophe and Bergerson, Joule and Brandt, Adam R and Masnadi, Mohammad S and Gordon, Deborah},
  journal={AGUFM},
  volume={2019},
  pages={A43M--2910},
  year={2019}
}

@article{liao2019batch,
  title={Batch-centric model for scheduling straight multisource pipelines},
  author={Liao, Q and Castro, PM and Liang, Y and Zhang, H},
  journal={AIChE Journal},
  volume={65},
  number={10},
  pages={e16712},
  year={2019},
  publisher={Wiley Online Library}
}

@article{caramanis1982optimal,
  title={Optimal spot pricing: Practice and theory},
  author={Caramanis, Michael C and Bohn, Roger E and Schweppe, Fred C},
  journal={IEEE Transactions on Power Apparatus and Systems},
  number={9},
  pages={3234--3245},
  year={1982},
  publisher={IEEE}
}

@book{schweppe2013spot,
  title={Spot pricing of electricity},
  author={Schweppe, Fred C and Caramanis, Michael C and Tabors, Richard D and Bohn, Roger E},
  year={2013},
  publisher={Springer Science \& Business Media}
}

@article{oneill2008towards,
  title={Towards a complete real-time electricity market design},
  author={O’Neill, Richard P and Fisher, Emily Bartholomew and Hobbs, Benjamin F and Baldick, Ross},
  journal={Journal of Regulatory Economics},
  volume={34},
  number={3},
  pages={220--250},
  year={2008},
  publisher={Springer}
}

@book{wood2013power,
  title={Power generation, operation, and control},
  author={Wood, Allen J and Wollenberg, Bruce F and Shebl{\'e}, Gerald B},
  year={2013},
  publisher={John Wiley \& Sons}
}

@article{litvinov2010design,
  title={Design and operation of the locational marginal prices-based electricity markets},
  author={Litvinov, Eugene},
  journal={IET generation, transmission \& distribution},
  volume={4},
  number={2},
  pages={315--323},
  year={2010},
  publisher={IET}
}

@article{lesieutre2005convexity,
  title={Convexity of the set of feasible injections and revenue adequacy in FTR markets},
  author={Lesieutre, BC and Hiskens, IA},
  journal={IEEE Transactions on Power Systems},
  volume={20},
  number={4},
  pages={1790--1798},
  year={2005},
  publisher={IEEE}
}

@article{oneill1979mathematical,
  title={A mathematical programming model for allocation of natural gas},
  author={O'Neill, Richard P and Williard, Mark and Wilkins, Bert and Pike, Ralph},
  journal={Operations Research},
  volume={27},
  number={5},
  pages={857--873},
  year={1979},
  publisher={INFORMS}
}

@inproceedings{rudkevich2017locational,
  title={Locational Marginal Pricing of Natural Gas subject to Engineering Constraints},
  author={Rudkevich, Alex and Zlotnik, Anatoly},
  booktitle={Proceedings of the 50th Hawaii International Conference on System Sciences},
  pages={3092--3101},
  year={2017}
}

@inproceedings{zlotnik2019optimal,
  title={Optimal Control for Scheduling and Pricing Intra-day Natural Gas Transport on Pipeline Networks},
  author={Zlotnik, Anatoly and Sundar, Kaarthik and Rudkevich, Aleksandr M and Beylin, Aleksandr and Li, Xindi},
  booktitle={2019 58th IEEE conference on decision and control (CDC)},
  pages={4877--4884},
  year={2019}
}

@inproceedings{rachford2000optimizing,
  title={Optimizing pipeline control in transient gas flow},
  author={Rachford Jr, Henry H and Carter, Richard G and others},
  booktitle={PSIG annual meeting},
  year={2000},
  organization={Pipeline Simulation Interest Group}
}

@article{zlotnik2016coordinated,
  title={Coordinated scheduling for interdependent electric power and natural gas infrastructures},
  author={Zlotnik, Anatoly and Roald, Line and Backhaus, Scott and Chertkov, Michael and Andersson, G{\"o}ran},
  journal={IEEE Transactions on Power Systems},
  volume={32},
  number={1},
  pages={600--610},
  year={2016},
  publisher={IEEE}
}

@article{wong1968optimization,
  title={Optimization of natural-gas pipeline systems via dynamic programming},
  author={Wong, Pl and Larson, R},
  journal={IEEE Transactions on Automatic Control},
  volume={13},
  number={5},
  pages={475--481},
  year={1968},
  publisher={IEEE}
}

@inproceedings{hajossy2014cooling,
  title={Cooling of a wire as the model for a rupture location},
  author={Hajossy, R and Mracka, I and Somora, P and Z{\'a}cik, T and others},
  booktitle={PSIG Annual Meeting},
  year={2014},
  organization={Pipeline Simulation Interest Group}
}

\end{document}